\documentclass{article}

\usepackage{amssymb,amsmath,mathrsfs}
\usepackage{pstricks,pst-node,pst-plot}

\def\mN{{\mathbb N}}
\def\cA{\mathcal A}

\def\cF{\mathcal F}
\def\cH{\mathcal H}

\def\cN{\mathcal N}

\def\ccG{\mathscr{G}}
\def\ccK{\mathscr{K}}

\def\ccP{\mathscr{P}}

\def\f{\varphi}
\def\s{\sigma}

\def\of{\widetilde{f}}

\def\G{\Gamma}
\def\DG{\Delta_{\ccG}}
\def\DP{\Delta_{\ccP}}
\def\FG{\cF_{\ccG}}
\def\FP{\cF_{\ccP}}

\def\homg{{\textsf{hom$_{\ccG}$}}}

\def\homm{{\textsf{Hom}}}

\def\sup{{\rm sup}\:}
\def\inf{{\rm inf}\:}

\def\proof{{\noindent \bf Proof :} }
\def\endproof{{\hfill $\square$} \vspace{2 mm}

}

\def\og{\leavevmode\raise.3ex\hbox{$\scriptscriptstyle\langle\!\langle$~}}
\def\fg{\leavevmode\raise.3ex\hbox{~$\!\scriptscriptstyle\,\rangle\!\rangle$}}

\newtheorem{theo}{Theorem}[section]
\newtheorem{cor}{Corollary}[section]
\newtheorem{de}{D\'efinition}[section] 
\newtheorem{lem}{Lemma}[section]    
\newtheorem{prop}{Proposition}[section] 
\newtheorem{rmk}{Remark}[section]

\newtheorem{ex}{Example}[section]

\def\redd{\searrow \!\!\stackrel{d}{}}
\def\redwd{\searrow \!\!\stackrel{wd}{}}
\def\expd{\stackrel{d}{} \!\!\nearrow }

\def\redred{\searrow\!\!\searrow\:}

\def\domd{\vdash\!\!\!\!^{^d}\:}

\def\fix{{\rm Fix}\:}
\def\sdr{strong deformation retract}

\def\ms{\medskip}
\def\ss{\smallskip}

\makeatletter
\def\timenow{\@tempcnta\time
  \@tempcntb\@tempcnta
  \divide\@tempcntb60
  \ifnum10>\@tempcntb0\fi\number\@tempcntb
  \multiply\@tempcntb60
  \advance\@tempcnta-\@tempcntb
  :\ifnum10>\@tempcnta0\fi\number\@tempcnta}
\makeatother

\usepackage{calc}
\newcounter{hours}\newcounter{minutes}

\begin{document}


\centerline{\Large Foldings in graphs and relations with simplicial complexes
and posets}
\vspace{8 mm}

\centerline{\large Etienne Fieux\footnote{ Institut de Math\'ematiques de Toulouse, 
Universit\'e Paul Sabatier, 118 Route de Narbonne, 
31062 Toulouse Cedex 09, France ; fieux@math.univ-toulouse.fr}, 
Jacqueline Lacaze\footnote{ Institut de Recherche en Informatique de Toulouse,
Universit\'e Paul Sabatier, 118 Route de Narbonne, 
31062 Toulouse Cedex 09, France ; jlacaze@irit.fr}}
\vspace{5 mm}

\centerline{Universit\'e Paul Sabatier, Toulouse, France}

\vspace{12 mm}



\begin{abstract}
We study dismantlability in graphs. In order to compare this notion to 
similar operations in posets (partially ordered sets) or in simplicial
complexes, we prove that a graph $G$ dismants on a subgraph $H$ if and only if
$H$ is a strong deformation retract  of $G$. Then, by looking at a triangle relating 
graphs, posets and simplicial complexes, we get a precise correspondence of the
various notions of dismantlability in each framework.
As an application, we study the link between the graph of morphisms
from a graph $G$ to a graph $H$ and the polyhedral complex Hom$(G,H)$;
this gives a more precise statement about 
 well known results concerning the polyhedral complex Hom$(G,H)$ and 
its relation with foldings in $G$ or $H$.
\vspace{3 mm}

\end{abstract}


\section{Introduction}

A vertex $g$ of a graph $G$ is said \textit{dismantlable}
 if there is another vertex $a$ in $G$ such that $N_G(x) \subset N_G(a)$  
 where $N_G(x):=\{y\in V(G), y\sim x\}$ is the open neighborhood of $x$.
This will be denoted $x \domd a$ 
and we will also say that $a$ dominates $x$.
The passage from $G$ to $G-x$ by deleting a dismantlable vertex $x$ is 
called a \textsl{folding} and denoted
$G \redd G-x$; the resultant graph 
$ G -x$ is called a \textsl{fold} of $G$.
A succession of foldings will be called a \textsl{dismantling}. 
If there is a dismantling
from a graph $G$ to a subgraph $H$, we say that
$G$ \textsl{is dismantlable on} and write
 $G \redd H$; this means that there is a \textsl{dismantling sequence} 
 $x_1,\ldots,x_k$ from $G$ to $H$, 
 i.e. $V(G) = V(H) \cup \{x_1,\ldots,x_k\}$ with $x_i$ dismantlable
in the subgraph induced by 
$V(H) \cup \{x_i,x_{i+1},\ldots,\ldots,x_k\}$ for $i=1,2,\ldots ,k$;
this will be also denoted $H \expd G$.
A reflexive graph $G$ is said \textsl{dismantlable} 
if it is dismantlable on a looped vertex.
Following \cite{hell-nesetril},
a graph whose every vertex is non dismantlable is called  {\sl stiff}.
 
It seems that the
 the first papers which focused on vertices whose open neighborhood
is included in the open neighborhood of another vertex\footnote{It is 
important to note that several papers 
(as \cite{anstee-farber},\cite{bfj08},\cite{ginsburg},\cite{larrion},\cite{quilliot83})
 take another definition of dismantlability : a
vertex $x$ is dismantlable if there is another vertex $a$ 
such that $N_G[x]\subset N_G[a]$ where $N_G[x]:=N_G(x)\cup \{x\}$
is the closed neighborhood of $x$. Of course, the two definitions
are the same when the graphs are reflexive.} are 
\cite{quilliot83} and  \cite{nowa-winkler}where it was proved
 independently that a reflexive graph is cop win if, and only if, 
 it is dismantlable. 
The reflexive bridged and connected graphs (a graph
is bridged if it  contains no isometric cycles of length greater than three;
in particular, the chordal graphs are bridged)
are  examples of dismantlable graphs (\cite{anstee-farber}).
 In this paper, the objective is to give a precise description of
 the relation between dismantlability in graphs and similar operations
 in partially ordered sets (posets) or in simplicial complexes.
In section 2, we give a characterization of foldings and dismantlings
by the way of morphisms and homotopies. 
The key result (Proposition \ref{Gdismant-and-sdr}) is
that a graph $G $ dismants on a subgraph  $H$ if, and only if,  
$H$ is a strong deformation retract of $G$. As a useful corollary, 
we get that if $G'$ and $G''$ are two subgraphs of a graphs $G$ such that
$G''$ is a subgraph of $G'$, $G\redd G'$ and $G\redd G''$ then
we can conclude that $G' \redd G''$ (Corollary \ref{dismant-subgraph}).

In the framework of posets, there is also  a very well known notion 
of dismantlability (most frequently named \textsl{irreducibility};
 see Section \ref{sectionGP}  for a brief discussion).  From 
the seminal paper \cite{stong}, we know that the dismantlings in posets
allow to describe the homotopy type of a poset (its \textsl{real}
homotopy type, i.e. the homotopy type of the poset considered 
as a topological space and not the homotopy type of its order
complex). Dismantlability in posets has been  studied in various
articles, in particular in relation with the fixed point property
(\cite{baclaw-bjorner},\cite{rival},\cite{schroder},\cite{walk84}).
It is known (\cite{bcf94},\cite{ginsburg}) that the dismantlability of 
a poset $P$ is equivalent to the dismantlability of its
\textsl{comparibility graph} (which will be called $Comp(P)$). 
In  \cite{ginsburg}, it was also  
proved that the dismantlability of a graph $G$
is equivalent to the dismantlability of 
the \textsl{poset of complete subgraphs} of $G$
(which will be called $C(G)$). 
In section 3, we will give a generalization of these results.
More recently (\cite{barmin09}), J. Barmak and A. Minian have introduced 
a notion of dismantlability in the category $\ccK$ of finite simplicial
complexes. We show in section 4 that 
it \og corresponds\fg~ to the dismantlability in
graphs under  natural functors relating $\ccG$ and $\ccK$. 

So, this gives a \og good\fg~ behaviour of a
 triangle $(\mathscr{G}^{\circ},\ccP,\ccK)$ in relation to the various 
 notions of dismantlability in $\ccG^{\circ}$, $\ccP$ or $\ccK$
 and, consequently, with the equivalences classes (named 
 \textsl{homotopy classes}) defined by the operation of dismantlability
 (section 5).
 A motivation for this question is given by the polyhedral complex
 $\homm(G,H)$ associated to two graphs $G$ and $H$.
 This construction is due to Lovasz after its pioneering work (\cite{lovasz})
where he solved the Kneser conjecture by using the simplicial complex 
$\cN(G)$, the \textsl{neighborhood complex} of $G$. Since the article 
\cite{babsonkozlov} (where the authors proved in particular that $\homm(K_2,G)$
and $\cN(G)$ have the same homotopy type), the $\homm$ complex
has became an important tool for determining lower bounds to the chromatic
number of certain graphs (see \cite{kozlov08} for a complete exposition
and more references). For obtaining topological information about
the polyhedral complex $\homm(G,H)$ (which is not in general a simplicial complex),
it is usual to look at its face poset $\FP(\homm(G,H))$ or at the order 
complex of its face poset, i.e. its barycentric subdivision
$Bd(\homm(G,H))=\DP(\FP(\homm(G,H)))$ (which is a simplicial complex).
On the other hand, the set of morphisms from $G$ to $H$
is the vertex set of a graph (called $\homg(G,H)$) and is also 
the vertex set of $\homm(G,H)$; we will study the relation 
 between the graph $\homg(G,H)$ and the polyhedral complex $\homm(G,H)$
 by using the triangle $(\mathscr{G}^{\circ},\ccP,\ccK)$ and
  regarding them in $\ccP$ (Proposition \ref{homg-homm-poset}).
 In particular, this gives another proof of a result describing the
 dismantlings on $\homm(G,H)$ induced by foldings on $G$ or $H$.
 However, this result which is usually formulated in terms of 
 simplicial complexes is formulated here in terms of graphs.
 \vspace{5 mm}

 \noindent {\bf Notations}
 In this paper, the graphs will be finite, undirected  and without
parallel edges. The vertex set of a graph $G$
is denoted $V(G)$. The set of these graphs will be denoted $\ccG$
and eventually considered as a category where a morphism 
$f : G \to G'$ from a graph $G$ to a graph $G'$ 
is an application from $V(G)$ to $V(G')$ which preserves adjacency 
($x\sim y \Longrightarrow f(x) \sim f(y)$);
$\ccG^{\circ}$ will denote
the subcategory obtained by retricting to reflexive graphs
(i.e., graphs $G$ such that $x \sim x$ for all $x$ in $V(G)$).

Let $G \in \ccG$. If $X$ is a subset of $V(G)$, the notation $G-X$
will indicate the subgraph of $G$ induced by the set of vertices
$V(G)\setminus X$.
In particular, if $x\in V(G)$,  $G-x$ will be
an abbreviated form of $G-\{x\}$
and $i_x : G-x \to G$ will denote the inclusion morphism.
If $x \domd a$, the folding $G \to G-x$ which sends $x$ to $a$
(and is the identity on $G-x$) will be denoted $r_{x,a}$.

The notation $G \expd G+y$ means that we have a added a vertex $y$ to $G$
in such a way that $y$ is dismantlable in the new graph.

\section{Morphisms}

In this section, we  characterize  foldings and dismantlings 
in terms of morphisms. Let $G,G'\in \ccG$. 
 The set  of morphisms from $G$ to $G'$ 
is the vertex set of a graph, denoted $\homg(G,G')$,
 where  $f\sim f'$ in $\homg(G,G')$ if and only if
$x\sim y$ in $G$ implies $f(x)\sim f'(y)$ in $G'$
(\cite{hhmnl},\cite{bcf94});
this graph is reflexive  because $f \sim f$ 
means precisely that $f$ is a morphism of graph.
By an abuse of notation, \og $f \in \homg(G,G')$\fg~
will mean that $f$ is a morphism from $G$ to $G'$
(in place of $f\in V(\homg(G,G'))$).

\begin{rmk}\label{rmk-compatible-with-circ}
Let $G,G',G'' \in \ccG$, $f,f' \in \homg(G,G')$ and
$h,h'\in \homg(G',G'')$.
If $f\sim f'$ and $h\sim h'$, then $h\circ f \sim h'\circ f'$
because $x\sim y$ in $G$ implies $f(x)\sim f'(y)$ in $G'$
(by $f\sim f'$) which implies $h\circ f(x)\sim h'\circ f'(y)$ in $G'$
(by $h \sim h'$).
\end{rmk}

\subsection{Foldings and retraction}

An important class of morphisms is given by retractions. 
A retraction of a graph $G$ to a subgraph $H$ of $G$ is a morphism
$r : G \to H$ such that $r(x)=x$ for all $x$ in $V(H)$. So, a morphism
$r:G\to G$ such that $r\circ r=r$ is a retraction of $G$ to $r(G)$.
The results of this paragraph are based on the following remarks:

\begin{rmk}\label{rmk-f-sim-1}
a. Let $f \in \homg(G,G)$. If $f \sim 1_G$ (where $1_G$
is the identity morphism on $G$), then every vertex $x$ of $G$ 
verifies either $f(x)=x$, or $x \domd f(x)$.

b. In particular, if $f:G\to G$ is a retraction such that $f\sim 1_G$,
then $f \redd f(G)$.
\end{rmk}

\begin{figure}[h]
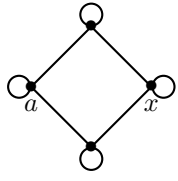

 \begin{center}
\psset{unit=0.8 cm}

\pspicture(-5,-1.5)(4,1.5)
\rput(-1,0){$\bullet$}
\pscircle(-1.2,0){0.2}
\uput[d](-1,0){$a$}
\rput(1,0){$\bullet$}
\pscircle(1.2,0){0.2}
\uput[d](1,0){$x$}
\rput(0,1){$\bullet$}
\pscircle(0,1.2){0.2}
\rput(0,-1){$\bullet$}
\pscircle(0,-1.2){0.2}
\pspolygon(-1,0)(0,1)(1,0)(0,-1)
\endpspicture
\caption{The retraction $G \to G-x$ (which sends $x$ to $a$) 
is not a folding}
\label{retraction-not-folding}
\end{center}
\end{figure}

We note that Remark \ref{rmk-f-sim-1}.a. implies that
$1_G$ is an isolated vertex in $\homg(G,G)$ when $G$ is a stiff graph
(this is a classical result used in \cite{bcf94}, \cite{docht}).
By definition, a folding is a retraction $G \to G-x$ which sends $x$ to a 
vertex $a$ which dominates $x$.
However, a general retraction $G \to G-x$ is not necessarily a folding 
(see Figure \ref{retraction-not-folding})

\noindent From  Remark \ref{rmk-f-sim-1}.a, 
we get the following characterization of foldings:

\begin{lem}\label{lemma-dismantlable-vertex}
Let $G \in \ccG$, $x\in V(G)$ and $f : G \to G-x$ a retraction; 
the following assertions are equivalent:
\vspace{2 mm}

1.  $f $ is a folding (i.e., $x \domd f(x)$)
\vspace{1 mm}

2. $i_x \circ f \sim 1_G$ (where $i_x$ 
is the inclusion $G-x \hookrightarrow G$)

\end{lem}


\noindent 
We conclude also from  Remark \ref{rmk-f-sim-1}.b that
foldings on  graphs induce dismantlability in graphs of morphisms:

\begin{prop} \label{foldings-on-homg}
1. If $x$ is dismantlable in $G$, then  $\homg(G,H) \redd \homg(G-x,H)$ 
(by identifying $\homg(G-x,H)$ with an induced subgraph of $\homg(G,H)$).

2. If $x$ is dismantlable in $H$, then  $\homg(G,H) \redd \homg(G,H-x)$ 
(by identifying  $\homg(G,H-x)$ with an induced subgraph of $\homg(G,H)$).
\end{prop}

\proof 1.  Let $x$ dismantlable in $G$ with $x\domd a$. Then, the map
$\Psi_{x,a} : \homg(G-x,H) \to \homg(G,H) $ 
defined by $\Psi_{x,a}(f)=f\circ r_{x,a}$ 
is an injective morphism of graphs
and we identify $\homg(G-x,H)$ with the subgraph 
$\Psi_{x,a}(\homg(G-x,H))$ of $\homg(G,H)$.
Let us denote $\Phi_{x} : \homg(G,H) \to \homg(G-x,H) $
the restriction morphism
defined by $\Phi_x(f)=f \circ i_x\equiv f_{\vert G-x}$. 
If $ f\in \homg(G-x,H)$, then $(\Phi_{x} \circ \Psi_{x,a})(f)=
(f\circ r_{x,a}) \circ i_x=f\circ (r_{x,a} \circ i_x)=f$;
so $\Phi_{x} \circ \Psi_{x,a}=1_{\homg(G-x,H)}$ and this means that
$\Psi_{x,a}\circ \Phi_{x} : \homg(G,H)\to \homg(G,H)$ is a retraction to
$\homg(G-x,H)$ identified with $\Psi_{x,a}(\homg(G-x,H))$. 
If $f \in \homg(G,H)$, $\Psi_{x,a}\circ \Phi_{x}(f)$
takes the same value as $f$ on vertices distinct from $x$ and 
takes the value $f(a)$
on $x$. Let $ f,f'\in \homg(G,H)$ with $f\sim f'$.
As  $i_x \circ r_{x,a} \sim 1_G$ 
(Lemma \ref{lemma-dismantlable-vertex}), we have
$f\circ i_x \circ r_{x,a} \sim f'$ by 
Remark \ref{rmk-compatible-with-circ} ; this proves that
$\Psi_{x,a}\circ \Phi_{x} \sim 1_{\homg(G,H)}$
and we conclude  $\homg(G,H) \redd \homg(G-x,H)$
by Remark \ref{rmk-f-sim-1}.b.

2. Similarly, if $x$ is dismantlable in $H$ with $x\domd b$, 
we denote $\Phi_{x,b} : \homg(G,H) \to \homg(G,H-x) $
the morphism of graphs defined by $\Phi_{x,b}(f)=r_{x,b}\circ f $ and
we identify  $\homg(G,H-x) $ with the induced subgraph
of  $\homg(G,H)$ given by its image under the injection
$\Psi_{x} : \homg(G,H-x) \to \homg(G,H) $
defined by $\Psi_{x}(f)=i_x \circ f$.  Then
$\Phi_{x,b} \circ \Psi_{x}=1_{\homg(G,H-x)}$ and
$\Psi_{x}\circ \Phi_{x,b} : \homg(G,H)\to \homg(G,H)$ is a retraction to
$\homg(G,H-x)$ identified with $\Psi_{x}(\homg(G,H-x))$. 
If $f \in \homg(G,H)$, $\Psi_{x}\circ \Phi_{x,b}(f)$
takes at a vertex $z$ the same value as $f$ when $f(z)\neq x$ 
and  the value $b$ when $f(z)=x$. It is easy to verify 
that $\Psi_{x}\circ \Phi_{x,b} \sim 1_{\homg(G,H)}$
and this proves  $\homg(G,H) \redd \homg(G,H-x)$.
\endproof

\subsection{Dismantlings and homotopy}

Morphisms give rise  to a notion of homotopy 
and it was noticed in
\cite{quilliot83} that a graph is dismantlable if and only if 
the identity morphism is homotopic to a constant morphism.
Following \cite{docht},
for $N \in \mN^*$,  $I_N$ 
is the reflexive graph with looped vertices $0,1,2,\ldots,N$
and adjacencies $0\sim 1\sim  2 \sim 3 \sim \ldots \sim N-1 \sim N$.

\begin{figure}[h]
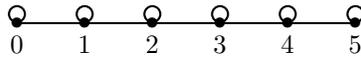

\begin{center}
\psset{unit=0.3}
\pspicture(10,-2)(5,2)
\psline(0,0)(15,0)
\rput(0,0){$\bullet$}\rput(3,0){$\bullet$}
\rput(6,0){$\bullet$}\rput(9,0){$\bullet$}\rput(12,0){$\bullet$}
\rput(15,0){$\bullet$}
\uput[d](0,0){0}\uput[d](3,0){1}
\uput[d](6,0){2}\uput[d](9,0){3}\uput[d](12,0){4}
\uput[d](15,0){5}
\pscircle(0,0.4){0.4}
\pscircle(3,0.4){0.4}
\pscircle(6,0.4){0.4}
\pscircle(9,0.4){0.4}
\pscircle(12,0.4){0.4}
\pscircle(15,0.4){0.4}
\endpspicture
\caption{The reflexive path $I_5$}
\label{reflexive-path}
\end{center}
\end{figure}

\noindent If $f,f' \in \homg(G,G')$, a homotopy from $f$ to $f'$ 
is a morphism of graphs $\cH :I_N \to \homg(G,G'),\: i \mapsto \cH_i$
such that $\cH_0=f$ and $\cH_N=f'$; 
this will be denoted $f\simeq f'$
and this means that $f$ and $f'$ are in 
the same connected component of $\homg(G,G')$.
 A subgraph 
 $G'$ of $G$ is a \textsl{\sdr} if there 
is a homotopy $\cH:I_N \to \homg(G,G)$ such that  $\cH_0=1_G$, 
${\cH_i}_{\vert G'}=1_{ G'}$ for all $i\in \{1,2,\ldots,N\}$ 
and $\cH_N:G\to G$ is actually a retraction to $G'$. 
 The following results will be useful in the sequel:

\begin{lem} \label{G2rdf-to-3rdf}
Let $G'' \subset G' \subset G$ inclusions of graphs.

1. If $G''$ is a \sdr ~ of $G'$ and $G'$ is a \sdr ~ of $G$, then
$G''$ is a \sdr ~ of $G$.

2. If  $G''$ is a strong deformation
retract of $G$ and $G'$ a retract of $G$, 
then $G''$ is a strong deformation retract of $G'$.
\end{lem}

\proof 1. Straightforward. 2.
Let $i_{G'} : G' \hookrightarrow G$ the inclusion and 
$r_{G'} :G \to G'$ a retraction of $G$  to $G'$.
If $\cH :I_N \to \homg(G,G)$ is a homotopy proving that 
$G''$ is a strong deformation retract of $G$, then
$\cH':I_N \to \homg(G',G')$
defined by  $\cH'_i =r_{G'} \circ \cH_i \circ i_{G'}$ is a homotopy 
proving that 
$G''$ is a strong deformation retract of $G'$.
\endproof

By Lemma \ref{lemma-dismantlable-vertex},
 a fold $G-x$ of $G$ is a \sdr ~of $G$; more generally,
 dismantlability is characterized by strong deformations:

\begin{prop}\label{Gdismant-and-sdr}
Let $H$ be a subgraph of a graph $G$.
Then, $G \redd H$ if, and only if,  
$H$ is a strong deformation retract of $G$.
\end{prop}

\proof
Let us suppose that $G \redd H$. 
This means that one can go from $G$ to $H$ by a composition
of foldings ; each fold being a \sdr, $H$ is a \sdr ~ of $G$ by  Lemma
\ref{G2rdf-to-3rdf} a.
If we suppose now that $H$ is a strong
deformation retract of $G$, we can use  an argument
similar to that used in the proof of Th\'eor\`eme 4.4 in \cite{bcf94}.
Let $\cH : I_N \to \homg(G,G)$ be a homotopy proving
that $H$ is a strong deformation retract of $G$.
If $H\neq G$, $\cH_N\neq 1_G$, 
and we can suppose $\cH_1 \neq 1_G$. 
If $a \in G$ is such that $\cH_1(a)\neq a$, 
then $a\notin H$, $a \domd \cH_1(a)$ in $G$ and $G \redd G-a$. 
For ${i \in\{0, 1,\ldots,N\}}$, define $\cH'_i : G-a \to G-a$ 
by $\cH'_i(x)=(r_{a,\cH_1(a)}\circ \cH_i)(x)$ for $x \in G-a$. 
Clearly, these morphisms define a homotopy 
$\cH' : I_N \to \homg(G-a,G-a)$,  $\cH'_1\sim 1_{G-a}$ and $\cH'_N$ is a 
retraction from $G-a$ to $H$. 
Thus, $H$ is a strong deformation retract of $G-a$.  
If $G-a \neq H$, we can iterate, and so $G \redd H$.
\endproof

\begin{cor}\label{dismant-subgraph}
Let $G'' \subset G' \subset G$ inclusions of graphs 
such that $G \redd G'$ and $G\redd G''$.
Then $G' \redd G''$.
\end{cor}

\proof Straightforward from Lemma \ref{G2rdf-to-3rdf} b. 
and Proposition \ref{Gdismant-and-sdr}.
\endproof

Let us recall that two graphs $G$ and $H$ 
are homotopically equivalent if there is $f \in \homg(G,H)$
and $g\in \homg(H,G)$ such that $g\circ f \simeq 1_G$ 
and $f\circ g \simeq 1_{H}$. In particular,
if $H$ is a \sdr~ of $G$,  $H$ and $G$ are homotopically equivalent. 
We mention the following well known result (two quite different proofs are
given in \cite{bcf94} and \cite{hell-nesetril}):

\begin{prop}\label{isom-two-stiff}
Let $G\in \ccG$ and $H,H'$ two stiff subgraphs 
such that $G\redd H$ and $G\redd H'$.
Then $H$ is isomorphic to $H'$.
\end{prop}

\proof
By Proposition \ref{G2rdf-to-3rdf},  $H$ and $H'$
are strong deformation retracts of $G$. So, $H$ and $H'$ are homotopically 
equivalent. Let $f\in \homg(H,H')$ and $g\in \homg(H',H)$ such that
$g\circ f \simeq 1_H$ and $f\circ g \simeq 1_{H'}$. As the graphs $H$ and $H'$ are stiff,
the connected components of $1_G$ and $1_H$ are reduced, respectively,
to $\{1_G\}$ and $\{1_H\}$. 
So, we conclude that $g\circ f = 1_H$ and $f\circ g = 1_{H'}$ and 
that $H$ and $H'$ are isomorphic.
\endproof

\section{Foldings versus $(\ccG,\ccP)$}\label{sectionGP}

Let $\ccP$ the category of finite posets. If $P,Q \in \ccP$, 
a morphism of posets $f :P \to Q$
is a map from $P$ to $Q$ which preserves the order (i.e. $x \leq y$ in $P$
implies $f(x)\leq f(y)$ in $Q$). An element
$p$ of  a poset $P$ will be called   \textit{dismantlable}
if either $P_{>p}:=\{y \in P, y>p\}$ has a least element 
or $P_{<p}:=\{y \in P, y<p\}$ has a greatest element. 
There are already various denominations
for this notion; the most classical are:   {\sl irreducible} 
(in many papers, following \cite{rival}), 
{\sl linear} and {\sl antilinear} (in \cite{stong}), 
{\sl upbeat points} and {\sl downbeat points} 
(\cite{may},\cite{barmin09});
in this paper, we adopt the denomination \og \textsl{dismantlable}\fg~
in order to emphasize the link with graphs.

Let $p$ a dismantlable point in $P$. If $a=\sup P_{<p}$
or $a=\inf P_{>p}$, $p$ will be said \textit{dominated by $a$}.
The deletion of the dismantlable element $x$, will be denoted 
$P\redd P\setminus \{x\}$ and $P \redd Q$ means that one can go from
the poset $P$ to a subposet $Q$ by successive deletions of
dismantlable elements.

\begin{prop}\label{Pdemont-sur-Fix}
Let $f :P\to P$ a morphism of posets map such that either $f\leq 1_P$ or $f\geq 1_P$. 
Then $P \redd \fix(f)$ where $\fix(f):=\{p\in P, f(p)=p\}$.
\end{prop}

\proof We suppose that $\fix(f)\neq P$ (i.e., $f\neq 1_P$) 
and we consider the case
 $f \leq 1_P$.
Let $x$ minimal in $P\setminus \fix(f)$. Let $y<x$
(for example, $y=f(x)$). Then $y=f(y)$ (by minimality of $x$ 
in $P\setminus \fix(f)$)
 and $y\leq f(x)$ (because $y<x\Rightarrow f(y)\leq f(x)$).
Thus, $f(x)$ is the greatest element of $P_{<x}$ and $x$ is dismantlable. 
So, we have $P \redd P\setminus \{x\}$.
Now, we define $\of : P\setminus \{x\} \to P\setminus \{x\}$ by $\of(y)=f(y)$ 
if $f(y) \neq x$ and $\of(y)=f(x)=f^2(y)$ if $f(y)=x$. 
Let us verify that $\of$ is a morphism of posets, i.e. 
$y \leq z \Longrightarrow \of(y) \leq \of(z)$, 
for all $y,z \in P\setminus \{x\}$. This is clear if either 
($f(y)\neq x$ and $f(z)\neq x$),
 or ($f(y)=f(z)=x$). If $f(y)=x$ and $f(z)\neq x$,
then $\of(y)=f^2(y) \leq f(y)\leq f(z)=\of(z)$ ($f^2(y) \leq f(y)$ follows from
$f\leq 1_P$ and $ f(y)\leq f(z)$ because $f$ is a morphism of posets).
Finally, if $f(y)\neq x$ and $f(z)= x$, we have $f(y) \leq f(z)$ 
(because $f$ is a morphism of posets),
so $f(y) < x$ (because $f(z)=x$ and $f(y)\neq x$). By minimality of 
$x$ in $P\setminus \fix(f)$, this means that $f(y) \in \fix(f)$. So we get
$f^2(y)=f(y)$ and $\of(y)=f(y)=f^2(y)\leq f^2(z)=f(x)=\of(z)$ 
(because $f^2$ is a morphism of posets).
It is clear that $\of : P\setminus \{x\} \to P\setminus \{x\}$ satisfies
$\of \leq 1_{P\setminus \{x\}}$ and that $\fix(\of)=\fix(f)$.
So, we can iterate the procedure and finally we get 
$P \redd P \setminus \fix(f)$.
The proof is similar if $f \geq 1_P$.
\endproof

\begin{rmk}
This proof is essentially the proof given 
by Kozlov in the particular case $f^2=f$
(\cite[Theorem 2.1]{kozlov06a} or \cite[Theorem 13.12]{kozlov08},
  where the conclusion is given in terms of simplicial complexes).
\end{rmk}

\subsection{Dismantlability and functor $Comp : \ccP \to \ccG^{\circ}$}

Let $P \in \ccP$. The comparability graph of $P$, denoted $Comp (P)$, 
is the graph
whose vertex set is $P$ with adjacencies $x \sim y$ if and only if
$x$ and $y$ are comparable in $P$. In particular, for every poset $P$,
$Comp(P)$ is a reflexive graph.
We will say that a graph $G$ is a cone with apex  $a$ if  $y \sim a$ 
for all $y\in V(G)$ (this definition implies that the apex is a looped vertex). 
The following facts are easy :

\begin{itemize}
\item If $x$ is a looped vertex of a graph $G$,  
then $x$ is dismantlable if, and only if,
$N_G(x)-x$ is a cone.
\item $Comp(P)-x=Comp(P\setminus \{x\})$.
\item $N_{Comp(P)}(x)-x=Comp(P_{>x}\cup P_{<x})$
\end{itemize}

\begin{prop}\label{dismant-PG}
Let $P,Q\in \ccP$.
If $P \redd Q$, then $Comp(P) \redd Comp(Q) $.
\end{prop}

\proof
Clearly, if $x$ is dominated by an element $a$ in $P$, then
 $x$ is dominated by the vertex  $a$ in $Comp(P)$. 
Consequently, $P \redd P\setminus \{x\}
 \Longrightarrow Comp(P) \redd Comp(P)-x=Comp(P\setminus \{x\})$
 and the proposition follows by iteration.
\endproof

\begin{figure}[h]
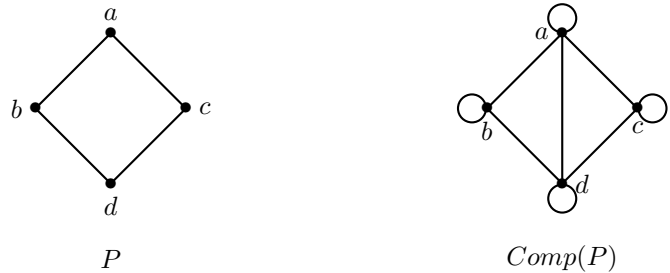

\begin{center}
  \pspicture(-3,-2)(4,2)
\rput(-1,0){$\bullet$}
\uput[l](-1,0){$b$}
\rput(1,0){$\bullet$}
\uput[r](1,0){$c$}
\rput(0,1){$\bullet$}
\uput[u](0,1){$a$}
\rput(0,-1){$\bullet$}
\uput[d](0,-1){$d$}
\pspolygon(-1,0)(0,1)(1,0)(0,-1)
\rput(0,-2){$P$}
\endpspicture
\pspicture(-2,-2)(4,2)
\rput(-1,0){$\bullet$}
\pscircle(-1.2,0){0.2}
\uput[d](-1,0){$b$}
\rput(1,0){$\bullet$}
\pscircle(1.2,0){0.2}
\uput[d](1,0){$c$}
\rput(0,1){$\bullet$}
\pscircle(0,1.2){0.2}
\uput[l](0,1){$a$}
\rput(0,-1){$\bullet$}
\pscircle(0,-1.2){0.2}
\uput[r](0,-1){$d$}
\pspolygon(-1,0)(0,1)(1,0)(0,-1)
\psline(0,1)(0,-1)
\rput(0,-2){$Comp(P)$}
\endpspicture
 \end{center}
\caption{$a$ and $d$ are non dismantlable in $P$ and dismantlable in $Comp(P)$}
\label{non-irreducible-in-P}
\end{figure}

Reciprocally, $Comp(P) \redd Comp(Q) $ doesn't imply in general 
$P \redd Q$ because a dismantlable vertex in $Comp(P)$ is not 
necessarily a dismantlable element in $P$ 
(see, for example,  the poset $P=\{a,b,c,d\}$ with $d<b,c<a$
given in  Figure \ref{non-irreducible-in-P}).

A poset will be called a \textsl{double cone with apex  $a$} if it admits 
an element $a$ comparable with all elements of the poset. It is clear that 
$P$ is a double cone if, and only if, $Comp(P)$ is a cone and this motivates the
following definition.

\begin{de}
An element $p$ of a poset $P$ is said \textit{weakly dominated by $a$} if
$P_{>p}\cup P_{<p}$ is a double cone with apex  $a$.
In this case,  $p$ will  be said \textit{weakly dismantlable}. We note 
$P\redwd P\setminus \{x\}$ the deletion of a weak dismantlable vertex 
and $P \redwd Q$ means that one can go from $P$ to a subposet $Q$ by 
successive deletions of weak dismantlable vertices.
\end{de}

In other terms, $p$ is weakly dominated by $a$ if $p$ and $a$ are comparable 
and if  every  element comparable with $p$ is also comparable with $a$. 
Of course, 
if $p$ is dominated by $a$, then  $p$ is weakly dominated by $a$
but the reverse implication is false in general (in the poset $P$ given in 
Figure \ref{non-irreducible-in-P},  
$d$ is weakly dominated by $a$ but is not dominated by $a$). 
Let $p$ an element of a poset $P$; the following assertions are equivalents :
\ms

1) $p$ is dominated by $a$ in $Comp(P)$\ss

2) $N_{Comp(P)}(p)-p$ is a cone with apex $a$\ss

3) $Comp(P_{>p}\cup P_{<p})$ is a cone with apex $a$ \ss

4) $P_{>p}\cup P_{<p}$ is a double cone with apex $a$\ss

5) $p$ is  weakly dominated by $a$  in $P$\ss
\ms

\noindent In particular, an element $p$ of a poset $P$ is weakly 
dismantlable in $P$ if, and only if, 
$p$ is dismantlable in $Comp(P)$ and the following equivalence is immediate :

\begin{theo}\label{theo-Comp}
Let $P,Q \in \ccP$. Then, 
$P \redwd Q ~\Longleftrightarrow~ Comp(P) \redd Comp(Q)$
\end{theo}

\subsection{Dismantlability and functor $C : \ccG \to \ccP$}

Let $G \in \ccG$. We recall that a  \textit{complete subgraph} $H$ of $G$
is an induced subraph of $G$ such that 
 $x\sim y$ for any distinct vertices $x$ and $y$
of $H$;   a complete subgraph of $G$ will be identified with its set
of vertices.  The poset of complete subgraphs of $G$, denoted $C(G)$,
is the poset given by the set of non empty complete subgraphs of $G$ 
with the  inclusion as order relation. 

\begin{theo}\label{theo-red-C}
Let $G \in \ccG$ and $H$ a subgraph of $G$ such that $G \redd H$.
If all vertices in $V(G)\setminus V(H)$ are looped, then 
$C(G) \redd C(H) $.

In particular, if $G \in \ccG^{\circ}$ and $H$ a subgraph of $G$, 
then $G \redd H \Longrightarrow C(G)\redd C(H)$.
\end{theo}

\proof Let $x\in V(G)\setminus V(H)$ be a  dismantlable and looped vertex
with $a$ which dominates $x$.
We define $f_1:C(G)\to C(G)$ by $f_1(c)=c \cup \{a\}$ if $x \in c$ 
and $f_1(c)=c$ if
$x\not\in c$; note that $f_1$ is well
defined because $x$ is  looped. Then $f_1 \geq 1_{C(G)}$  and, by 
Proposition \ref{Pdemont-sur-Fix}, $C(G)\redd {\rm Im}(f_1)$.
Now, let $f_2:{\rm Im}(f_1)\to {\rm Im}(f_1)$ 
 defined by $f_2(c)=c \setminus \{x\}$ if $x \in c$ 
and $f_2(c)=c$ if $x\not\in c$.
Then $f_2 \leq 1_{{\rm Im}(f_1)}$  and, by 
Proposition \ref{Pdemont-sur-Fix}, ${\rm Im}(f_1)\redd {\rm Im}(f_2)=C(G-x)$.
So, $C(G)\redd C(G-x)$ and the proposition follows by iterating the process.
\endproof

Now, before studying the reciprocal of Theorem \ref{theo-red-C}, we recall
that an element $p$ of a poset $P$ is an atom if $P_{<p}=\emptyset$ ;
the set of atoms of a poset $P$ will be denoted $\cA(P)$.
We introduce the 
applications
$$RUB :\ccP \to \ccG^{\circ}~~~~~~~~{\rm and}~~~~~~~~m : \ccP \to \ccG^{\circ}$$

\begin{itemize}
\item $RUB(P)$ is the \textit{reflexive upper bound graph} of $P$ :
$V(RUB(P))=P$ et and $p\sim q$ in $RUB(P)$ if there is a $z\in P$ 
such that $z\geq p$ and $z \geq q$ (in other words,
$p \sim q \Longleftrightarrow P_{\geq p,q}:=
P_{\geq p} \cap P_{\geq q} \neq \emptyset$).

\item $m(P)$ is the subgraph of $RUB(P)$ induced by $\cA(P)$, the set
of atoms  of $P$ (i.e., $V(m(P))=\cA(P)$ 
and $a\sim b \in m(P)$ if there is a $p\in P$ 
such that $p\geq a$ and $p\geq b$).
\end{itemize}

\begin{prop} \label{RUB-dd-to-m}
For all $P\in \ccP$, $RUB(P) \redd m(P)$.
\end{prop}

\proof As $m(P)$ is the subgraph of $RUB(P)$ induced by the set
 of atoms $\cA(P)$, it suffices to prove that every 
 vertex in $V(RUB(P))\setminus V(m(P))$
(i.e., every element of $P$ which is not an atom) 
is dominated by a vertex of $m(P)$.
Let $q \in V(RUB(P))\setminus V(m(P))=P \setminus \cA(P)$. It is immediate 
that $q \domd x$ for every  vertex 
$x \in V(m(P))=\cA(P)$ such that $x<q$
(because $z \sim q \Longleftrightarrow P_{\geq z,q}\neq \emptyset
\Longrightarrow P_{\geq z,x}\neq \emptyset \Longrightarrow z \sim x$).
\endproof

\begin{prop} \label{RUB-et-dd}
Let $P$ in $\ccP$ and $x$ dismantlable in $P$.
Then, $RUB(P) \redd RUB(P\setminus \{x\})$.

As a consequence, $P \redd Q \Longrightarrow RUB(P)\redd RUB(Q)$.
\end{prop}

\proof Let us suppose that $x$ is dominated by $a$ in $P$.
First, we verify that $x$ is dominated by $a$ in $RUB(P)$.
So, let $y \in P$ such that 
$y \sim x$ in $RUB(P)$. If $y=x$, then $x\sim a$ in $RUB(P)$
because $P_{\geq y,a}=P_{\geq x,a}\neq \emptyset$.
If $y\neq x$ and $z\in P_{\geq y,x}$, 
then $z\in P_{\geq y,a}$ 
(because $z\geq x$ and $x$ is dominated by $a$) ; so, $y \sim a$
and $x$ is also dominated by $a$ in  $RUB(P)$. 
Hence we have $RUB(P) \redd RUB(P)-x$. 
Now, we compare the graphs $RUB(P)-x$ and $RUB(P\setminus \{x\})$. 
They have the same vertex sets and 
clearly $RUB(P\setminus \{x\})$ is a subgraph of $RUB(P)-x$ 
(if $P_{\geq y,z}\neq \emptyset $ in $P\setminus \{x\}$,
we have also $P_{\geq y,z}\neq \emptyset $ in $P$). 
Now, let us suppose that $y\sim z$ in $RUB(P)-x$; this means that
$P_{\geq y,z} \neq \emptyset$. If $x$ is in $P_{\geq y,z}$,
then $a$ is also in $P_{\geq y,z}$; so, 
$P_{\geq y,z} \cap (P\setminus \{x\}) \neq \emptyset$
and this proves that 
$y\sim z$ in $RUB(P\setminus \{x\})$. 
In conclusion, $RUB(P)-x= RUB(P\setminus \{x\})$ and 
$RUB(P) \redd RUB(P\setminus \{x\})$.
\endproof

 Let us denote by $G^{\circ}$ the reflexive graph obtained from a graph $G$
by adding loops to its non looped vertices.
We note that, by identifying  $\cA(C(G))$ with $V(G)$, we get $m(C(G))=G^{\circ}$
for every $G \in \ccG$.

\begin{theo}\label{theo-C-red}
Let $G\in \ccG$ and $H$  a  subgraph of $G$ such that $C(G)\redd C(H)$.
 Then $G^{\circ} \redd H^{\circ}$.
 
 In particular, if $G \in \ccG^{\circ}$ and $H$ is a subgraph of $G$, 
then $C(G)\redd C(H) \Longrightarrow G \redd H$.
\end{theo}

\proof 
By Proposition \ref{RUB-dd-to-m}, we have two dismantling
$f_G : RUB(C(G))\redd m(C(G))= G^{\circ}$
and $f_H: RUB(C(H))\redd m(C(H))=H^{\circ}$.
There is also a dismantling $\f : RUB(C(G)) \redd RUB(C(H))$
from $C(G)\redd C(H)$ and Proposition \ref{RUB-et-dd}.
So, we have the following diagram:
$$
\begin{psmatrix}[colsep=2cm,rowsep=2cm]
RUB(C(G)) & m(C(G))=G^{\circ}\\
RUB(C(H)) & m(C(H))=H^{\circ}
\psset{arrows=->,labelsep=5pt,nodesep= 0.5cm}
\ncline{1,1}{1,2}^{f_G}
\ncline{1,1}{2,1}<{\f}
\ncline{2,1}{2,2}^{f_H }
\ncline[linestyle=dashed]{1,2}{2,2}
\end{psmatrix}
$$
The conclusion  $G^{\circ} \redd H^{\circ}$ 
follows from Corollary \ref{dismant-subgraph} applied 
to $(G,G',G'')=(RUB(C(G)),G^{\circ},H^{\circ})$.
\endproof

\section{Foldings versus $(\ccG,\ccK)$}

\subsection{Dismantlability in $\ccK$}

Let $K$ be the category of finite simplicial complexes (cf. \cite{kozlov08} 
for a  reference textbook) and let $K \in \ccK$. 
If $\sigma$ is a simplex of $K$, we write $\sigma \in K$.
A simplicial complex $K$ is a \textit{simplicial cone} if there is a 
subcomplex $L$ and a vertex $a$ of $K\setminus L$
such that the set of simplices of $K$ is
$\{\{a\},\sigma,\{a\} \cup \s, \s \in L\}$; 
in this case, $K$ is denoted $aL$.
Let us recall the following definitions for a vertex $x$ of $K$:

\begin{itemize}
\item ${\rm star}^o_K(x):=\{\sigma \in K,\: x\in  \sigma \}$
\item ${\rm lk}_K(x):=\{\sigma \in K,\: \{x\} \cup \sigma \in K~{\rm and}~
x \notin \sigma\}$
\item ${\rm star}_K(x):=\{\sigma \in K,\: 
\{x\} \cup \sigma \in K\}={\rm star}^o_K(x) \cup {\rm lk}_K(x)$
\item $K-x:=\{\sigma \in K,\: x\notin  \sigma \}$.
\end{itemize}

We note that a simplicial complex $K$ is a simplicial cone if, and only if,
one can write $K=xL$ with $L=K-x$ for some vertex $x$.
In \cite{barmin09}, a notion of dismantlability is defined in the framework 
of simplicial complexes.
A vertex $x$ of a simplicial complex $K$ is said \textit{dominated} by the vertex $a$
of $K$ if ${\rm lk}_K(x)$ is
a simplicial cone $aK'$ for some subcomplex $K'$ of $K$; 
in this case,  the deletion of the vertex $x$ in $K$ is called
\textit{an elementary strong collapse} and denoted $K \redred K-x$. 
A strong collapse, denoted  $K \redred L$,
is the succession of elementary strong collapses.
In this paper, by analogy with the situation in graphs and posets,
a dominated vertex in a simplicial complex $K$ will be said
\textit{dismantlable} in $K$.

\begin{rmk}
In \cite{cy}, the authors introduce the notion of \textsl{linear coloring} on simplicial
complexes. 
The Theorem 6.2 of \cite{cy} shows
that the notion of LC-reduction in \cite[\S 6]{cy} and the notion of strong reduction 
defined in \cite{barmin09} are equivalent.




\end{rmk}

\subsection{Dismantlability and functor $\DG : \ccG \to \ccK$}

Let $G \in \ccG$. We recall that $\DG(G)$ (sometimes called \textit{the clique
complex} of $G$) is the simplicial complex whose 
simplices are given by sets of vertices of complete subgraphs of $G$.
The following facts are easy :

\begin{itemize}
\item  If $G$ is a reflexive graph, then $G$ is a cone if, and only if, 
$\DG(G)$ is a simplicial cone.
\ss

\item  For every vertex $x$ of a graph $G$, $\DG(N_G(x)-x)={\rm lk}_{\DG(G)}(x)$.
\end{itemize}

\begin{lem}\label{dismantlable-vertex-G-K}
Let $G \in \ccG$, $a,x\in V(G)$ such that $a\neq x$ and $x$ looped.
Then, $x$ is dominated by $a$ in $G$
if, and only if, $x$ is dominated by $a$ in $\DG(G)$.
\end{lem}

\proof
Let $x$ a looped vertex. If $x\domd a$, then $N_G(x)-x$ is a cone with apex $a$
and ${\rm lk}_{\DG(G)}(x)=\DG(N_G(x)-x)=
aL$ with $L=\bigl(\DG(N_G(x)-x)\bigr)-a$, i.e. $x$ is dominated
by $a$ in $\DG(G)$.
Conversely, if ${\rm lk}_{\DG(G)}(x)=\DG(N_G(x)-x)$ is a simplicial cone $aL$,
then necessarily  $y\sim a$ for all $y \in N_G(x)-x$ and $x\sim a$;
in other terms, $N_G(x)\subset N_G(a)$, i.e.  $x\domd a$.
\endproof

\begin{theo}\label{theo-Delta}
Let $G ,H\in \ccG^{\circ}$. Then,
$G \redd H ~\Longleftrightarrow~\DG(G) \redred \DG(H) $.
\end{theo}

\proof Follows by iteration of Lemma 
\ref{dismantlable-vertex-G-K}.\endproof

\subsection{Dismantlability and functor $\FG : \ccK \to \ccG^{\circ}$}

Let $K$  a simplicial complex. The \textsl{face graph} (\cite{bfj08})
$\FG(K)$ of $K$
is the reflexive graph whose vertices
are the non empty simplices of $K$
with an edge between two simplices if one contains the other.
If $x$ is a vertex of $K$, $\{x\}$ will denote the same vertex as
a 0-simplex of $K$ or as a vertex of $\FG(K)$.
More generally, if $\s$ is a simplex of $K$, we also denote $\s$ 
the corresponding vertex of $\FG(K)$.

\begin{theo}\label{theo-red-Gamma}
Let $K,L\in \ccK$. Then,
$ K \redred L ~\Longrightarrow~\FG(K) \redd \FG(L) $.
\end{theo}

\proof
It is sufficient to prove  
$K \redred K-x ~\Longrightarrow~\FG(K)\redd \FG(K-x)$.
 As $V(\FG(K))\setminus V(\FG(K-x))={\rm star}^o_K(x)$,
we have to verify that one can dismant, one by one, all the elements of
${\rm star}^o_K(x)$ when ${\rm lk}_K(x)$ is a cone. 
 So, let $x$ a dismantlable vertex in $K$ and $a$ a vertex
which dominates $x$ in $K$; we have
$ {\rm star}^o_K(x) =\G_x \cup \G_{x,a}$ with
 $\G _{x} :=\{\s \in K, x \in \s ~{\rm and}~ a \not\in \s\}$ and 
$\G_{x,a} :=\{\s \in K, x \in \s ~{\rm and}~ a \in \s\}$. 
As the neighborhood in $\FG(K)$ of a  maximal simplex 
$\s$ of $\G_x$ is $\{\{a\}\cup \sigma\}\cup \{ \tau, \tau \subset \s\}$,
we have $\s \domd \{a\}\cup \sigma$ in $\FG(K)$.
 So, all maximal simplices of $\G_x$ are dismantlable
and, when they have been deleted,
the maximal simplices of the resulting subset of $\G_x$ are also dismantlable
by the same argument and the iteration of this procedure
showes that all vertices of $\G_{x}$ are dismantlable (the procedure ends
when the 0-simplex $\{x\}$ is dominated by the 1-simplex $\{a,x\}$).  
Next, it remains to prove that one can dismant all vertices of $\G_{x,a}$. 
This follows  from the existence of a similar procedure
to the precedent, in the reverse order. 
First, the vertex $\{x,a\}$ is dominated by $a$.
Next, after the removing of $\{x,a\}$, 
vertices of type $\{a,x,y\}$ are dominated by $\{a,y\}$ and
after the removing of these vertices, 
vertices of type $\{a,x,y,z\}$ are dominated by $\{a,y,z\}$
and so on, until all vertices of $\G_{x,a}$ have been deleted.
\endproof

\begin{rmk}
There is an obvious morphism $f \circ g: \FG(K)\to \FG(K-x)$ 
where $g : \FG(K) \to \FG(K)- \G_x $ 
is defined by $g(\s)=\{a\}\cup \s$ on $\G_x$ 
and $g(\s)=\s$ otherwise and $f : \FG(K)- \G_x \to \bigl( \FG(K) 
-\G_x\bigr) -\G_{x,a}=\FG(K-x) $ is defined by 
$f(\s)= \s \setminus \{x\}$ on $\G_{x,a} $  and $f(\s)=\s$ otherwise.
 Nevertheless,
in general, we don't have $g \sim 1_{\FG(K)}$, 
nor $f\sim 1_{ \FG(K)- \G_x }$ 
and the preceding proof shows the necessity
of deleting the vertices of ${\rm star}^o_K(x) $  in a certain order.
\end{rmk} 

To establish the reciprocal statement of Theorem \ref{theo-red-Gamma}, 
we need two lemmas.

\begin{lem}\label{lemme_preparatoire}
Let $K\in \ccK$ and $L$  a subcomplex of $K$ such that $\FG(K)\redd \FG(L)$.
If $\sigma$  is a maximal simplex of $K$ which appears in a dismantling 
sequence from $\FG(K)$ to $\FG(L)$, then there is a $0$-simplex 
$\{x\}$ with $x\in \s$ which appears  before $\s$ in the same dismantling sequence.
\end{lem}

\proof 
Let us suppose that $\s$ is a maximal simplex of $K$
 which appears in a dismantling sequence from $\FG(K)$ to $\FG(L)$.
This means that after having removed some vertices, we get a subgraph
$\cF'$ of $\FG(K)$ and there is a simplex $\s'$ which dominates $\s$
in $\cF'$.
As $\s$ is a maximal simplex and $\s \sim \s'$, we must have 
$\s'\varsubsetneq \s$.
Now, let $x\in \s$ ; $\s \domd \s'$ implies $\{x\} \sim \s'$,
i.e. $x \in \s'$.
In particular, if no vertex of $\sigma$ has been dismantled,  then
$\s \subset \s'$.  But this contradict $\s'\varsubsetneq \s$. 
So, there must be at least one vertex of $\sigma$
which has been dismantled before $\s$.
\endproof

\begin{lem}\label{dismant-vertex}
Let $K\in \ccK$ and $L$  a subcomplex of $K$ such that $\FG(K)\redd \FG(L)$.
If $\{x\}$ is the first 0-simplex  dismantled in a  dismantling sequence
from $\FG(K)$ to $\FG(L)$, then $x$ is dismantlable in $K$.
\end{lem}

\proof
Let $\{x\}$  be the first 0-simplex  dismantled in a  dismantling sequence
from $\FG(K)$ to $\FG(L)$ 
and $\sigma$ a simplex such that $\{x\} \domd \sigma$
in the dismantling process.
We will show that every element  of $\sigma$ dominates $x$ in $K$.
So, let us take $a \in \sigma$, $a\neq x$ and $\tau  \in {\rm lk}_K(x)$.
We have to prove that $\tau \cup \{a\}$ is a simplex of ${\rm lk}_K(x)$.
Let $\tau_{\max}$ be a maximal simplex of $K$ containing $\tau \cup \{x\}$;
by Lemma \ref{lemme_preparatoire}, we know that $\tau_{\max}$ has not
been dismantled before $x$.  As $x \in \tau_{\max}$ 
and $\{x\}\domd \s$, 
we conclude that $\s$ is adjacent to $\tau_{\max}$, i.e. 
$\s \subset  \tau_{\max}$ (because $\tau_{\max}$ is maximal).
Consequently, $a \in \tau_{\max}$ and
$\tau \cup \{a,x\} \subset \tau_{\max}$; this shows
that $\tau \cup \{a\}$ is a simplex of ${\rm lk}_K(x)$.
\endproof

\begin{theo} \label{theo-Gamma-red}
Let $K\in \ccK$ and $L$  a subcomplex of $K$ such that $\FG(K)\redd \FG(L)$.
 Then $K \redred L$.
\end{theo}

\proof 
By Lemma \ref{dismant-vertex}, 
we know that there exists a vertex  $x$ of $ K-L$
such that $K \redred K-x$. Now, from Theorem \ref{theo-red-Gamma}, we get
a dismantling $f_x:\FG(K)\redd \FG(K-x)$. So, with the hypothesis
of a dismantling $\f:\FG(K) \to \FG(L)$, we have the following triangle:
$$
\begin{psmatrix}[colsep=2cm,rowsep=2cm]
\FG(K) &\FG(K-x)\\
 & \FG(L)
\psset{arrows=->,labelsep=5pt,nodesep= 0.5cm}
\ncline{1,1}{1,2}^{f_x}
\ncline{1,1}{2,2}<{\f}
\ncline[linestyle=dashed]{1,2}{2,2}
\end{psmatrix}
$$
which allows to  conclude $\FG(K-x) \redd \FG(L)$ from Corollary 
\ref{dismant-subgraph} because $\FG(L)$ is a subgraph of $\FG(K-x)$.
Now, we iterate the argument with $\FG(K-x)$. 
The iteration ends 
when all 0-simplices  which are not vertices of $\FG(L)$ have been dismantled
and this proves that $K \redred L$.
\endproof

\begin{rmk}
We also deduce from the proof of Theorem \ref{theo-Gamma-red} 
that a dismantling sequence from $K$ to $L$ is obtained by keeping 
the 0-simplices (or vertices of $K$) 
 in a dismantling sequence from $\FG(K)$ to $\FG(L)$. 
\end{rmk}

\section{Homotopy classes and the triangle $(\mathscr{G}^{\circ},\ccP,\ccK)$ }

\subsection{Posets, simplicial complexes and dismantlability}

The order complex of a poset $P \in \ccP$
is the simplicial complex $\DP(P)$ 
whose simplices are given by the chains of $P$.
First, we note the elementary facts:

\begin{itemize}
\item $\Delta_{\ccP}(P_{>x} \cup P_{<x})={\rm lk}_{\DP(P)}(x)$.
\item A poset $P$ is a double cone with apex $a$ (i.e., 
$P=P_{>a}\cup P_{<a} \cup \{a\})$ if, and only if, 
$\Delta_{\ccP}(P)$ is a simplicial cone with apex $a$.
\end{itemize}

\noindent As a consequence of these facts, an element $x$ of $P$ 
is weakly dismantlable if, and only if, 
$x$ is dismantlable in $\Delta_{\ccP}(P)$ and:

\begin{theo}\label{theo-DP} 
Let $P,Q \in \ccP$.
Then, $P \redwd Q ~\Longleftrightarrow~\Delta_{\ccP}(P)
 \redred\Delta_{\ccP}(Q) $.
\end{theo}

\begin{rmk}
We know from \cite[Theorem 4.14.a]{barmin09} that
$P \redd Q $ implies $\Delta_{\ccP}(P) \redred\Delta_{\ccP}(Q) $;
the  example of the poset $P$  given in 
Figure \ref{non-irreducible-in-P} 
($d$ is dominated by $a$ in $\DP(P)$  but not dominated in $P$) shows that 
the reciprocal statement is not true in general.
\end{rmk}

Let $K$ a simplicial complex.
The \textsl{face poset} $\FP(K)$ of $K$ is the poset given by the set
of non empty simplices of $K$ with the inclusion as order relation. 
From \cite[Theorem 4.14.b]{barmin09}, 
we know that
$K \redred L ~\Longrightarrow~\FP(K) \redd \FP(L) $;
the reciprocal statement is true:

\begin{theo}\label{theo-FP} 
Let $K,L\in \ccK$. If 
 $\FP(K) \redd \FP(L)$, then $K \redred L $.
\end{theo}

\proof 
Let us suppose that  $\FP(K) \redd \FP(L)$ in $\ccP$. 
By Proposition \ref{dismant-PG} and identity $Comp \circ \FP=\FG$,
$\FG(K) \redd \FG(L)$ in $\ccG^{\circ}$ and, by  Theorem \ref{theo-Gamma-red}, 
$K \redred L$.
\endproof

\subsection{Homotopy classes}

Addition or deletion of dismantlable vertices define 
an equivalence relation in $\ccG$: $[G]_d=[H]_d$
if there is  in $\ccG$ a sequence 
$G=J_0,J_1,J_2,\ldots,J_{n-1},J_n=H$ such that 
$J_i\redd J_{i+1}$ or $J_i \expd J_{i+1}$ or 
$J_i\cong J_{i+1}$ ($J_i$ and $J_{i+1}$ 
are isomorphic graphs) for $i=0,1,\ldots,N-1$. 
The equivalence class $[G]_d$ of $G$ will be called the 
{\sl d-homotopy type} of $G$.

The term \textsl{homotopy} is given here by analogy with the 
equivalence class $[P]_d$ of a  poset $P$ in $\ccP$
which is defined in a similar way by dismantlings in $\ccP$.
It is well known (\cite{stong})
 that $[P]_d$ is actually the homotopy class of the poset $P$
considered as a topological space (with $\{P_{\leq x},x \in P\}$  
as a base of neighborhoods). 
In a similar way and following \cite[Definition 2.1]{barmin09},
two simplicial complexes $K$ and $L$
have the same \textsl{strong homotopy type}  
if one can go from $K$ to $L$ by a succession of strong collapses
or strong expansions.

\begin{rmk}\label{rigidity-d-homotopy}
The $d$-homotopy type is quite rigid. The example (see Figure \ref{the_Cn}) of the 
reflexive cycles $C_n^{\circ}$ shows the important gap with the $s$-homotopy
(two graphs $G$ and $H$ have the same $s$-homotopy type if, and only if,
the simplicial complexes $\DG(G)$ and $\DG(H)$ have the same 
simple homotopy type, cf. \cite{bfj08}).

\begin{figure}[h]
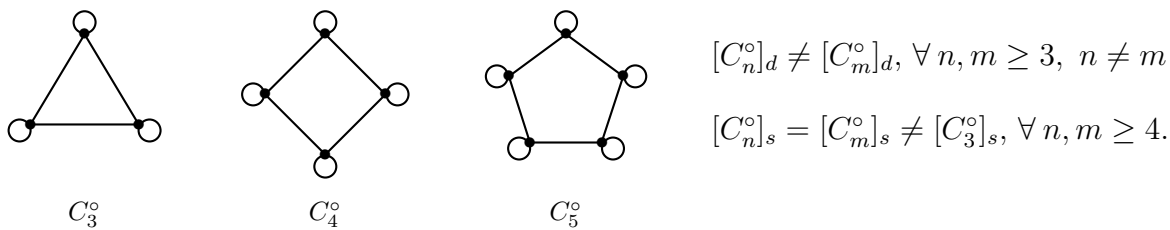

\begin{center}
\psset{unit=0.8cm}
\pspicture(6,-2)(6,2)
\pspolygon(0,1)(0.9,-0.5)(-0.9,-0.5)
\rput(0,1){$\bullet$}
\pscircle(0,1.2){0.2}
\rput(-0.9,-0.5){$\bullet$}
\pscircle(-1.08,-0.6){0.2}
\rput(0.9,-0.5){$\bullet$}
\pscircle(1.08,-0.6){0.2}
\rput(0,-2){$C_3^{\circ}$}
\endpspicture
\pspicture(5,-2)(2,2)
\pspolygon(-1,0)(0,1)(1,0)(0,-1)
\rput(-1,0){$\bullet$}
\pscircle(-1.2,0){0.2}
\rput(1,0){$\bullet$}
\pscircle(1.2,0){0.2}
\rput(0,1){$\bullet$}
\pscircle(0,1.2){0.2}
\rput(0,-1){$\bullet$}
\pscircle(0,-1.2){0.2}
\rput(0,-2){$C_4^{\circ}$}
\endpspicture
\pspicture(1,-2)(2,2)
\pspolygon(0,1)(-0.95,0.3)(-0.6,-0.8)(0.59,-0.8)(0.95,0.3)
\rput(0,1){$\bullet$}
\pscircle(0,1.2){0.2}
\rput(-0.95,0.3){$\bullet$}
\pscircle(-1.15,0.3){0.2}
\rput(-0.6,-0.8){$\bullet$}
\pscircle(-0.78,-0.9){0.2}
\rput(0.6,-0.8){$\bullet$}
\pscircle(0.78,-0.9){0.2}
\rput(0.95,0.3){$\bullet$}
\pscircle(1.15,0.3){0.2}
\rput(0,-2){$C_5^{\circ}$}
\rput(6.3,0){\large
\begin{tabular}{l}
$[C_n^{\circ}]_d\neq [C_m^{\circ}]_d$, $\forall \: n,m\geq 3,~n\neq m$\\
~~~\\
$[C_n^{\circ}]_s=[C_m^{\circ}]_s\neq [C_3^{\circ}]_s$, 
$\forall \: n,m\geq 4$.
\end{tabular}
}
\endpspicture
\end{center}
\caption{$d$-homotopy classes and $s$-homotopy classes}
\label{the_Cn}
\end{figure}

\end{rmk}

\begin{prop}\label{wd=d}
Let $P \in \ccP$ and $x$ a weak dismantlable element in  $P$.
Then $[P]_d=[P\setminus \{x\}]_d$.
\end{prop}

\proof 
As $x$ is a weak dismantlable element in  $P$, it exists
an element $a$ comparable with all elements of $P_{>x}\cup P_{<x}$.
Let us suppose that $a>p$ and let $y\in \max [x,a[$
where $[x,a[=P_{\geq x}\cap P_{<a}$. It is easy to see that $y\domd a$
(indeed, if $z>y$, then $z>x$, so $z$ is comparable with $a$;
but $z<a$ would contradict  $y\in \max [x,a[$, so $z\geq a$).
So, one can remove by dismantlability all maximal elements of 
$[x,a[$ and the iteration of this reasoning until $x$ is removed
proves that $P \redd Q$ and $P\setminus \{x\} \redd Q$ with $Q=P\setminus [x,a[
=(P\setminus \{x\})\setminus ]x,a[$. If we suppose that $a<p$,
we get $P \redd Q'$ and $P\setminus \{x\} \redd Q'$ with $Q'=P\setminus ]a,x]
=(P\setminus \{x\})\setminus ]a,x[$.
\endproof

 \begin{rmk}
As a useful consequence of Proposition \ref{wd=d}, if two posets
$P$ and $Q$ are such that $P \redwd Q$, then $[P]_d=[Q]_d$. In other terms,
the weak dismantlability preserves the homotopy type in $\ccP$.
\end{rmk}

\subsection{The triangle $(\mathscr{G}^{\circ},\ccP,\ccK)$}

\begin{figure}[h]
 \begin{center}
\psset{unit=0.6 cm}
\pspicture(-2,-0.5)(6,6.5)
\rput(0,0){\large $\mathscr{G}^{\circ}$}
\rput(0,6){\large $\ccP$}
\rput(5,3){\large $\ccK$}
\psline{->}(1,5.1)(4,3.6)
\rput(2.5,5.5){\large $\FP$}
\psline{->}(4,4.2)(1,5.7)
\rput(2.5,3.8){\large $\Delta_{\ccP}$}
\psline{->}(-0.3,4.8)(-0.3,1.2)
\rput(-1.4,3){\large $Comp$}
\psline{->}(0.3,1.2)(0.3,4.8)
\rput(1.1,3){\large $C$}
\psline{->}(1,0.9)(4,2.4)
\rput(2.5,2.3){\large $\Delta_{\mathscr{G}}$}
\psline{->}(4,1.8)(1,0.3)
\rput(2.5,0.5){\large $\FG$}
\endpspicture
\caption{The triangle $(\mathscr{G}^{\circ},\ccP,\ccK)$}
\label{GPK}
\end{center}
\end{figure}

\noindent The  functors in the triangle  $(\mathscr{G}^{\circ},\ccP,\ccK)$
are compatible with the various homotopy classifications
(d-homotopy type in $\ccG$, homotopy type in $\ccP$ and strong homotopy
type in $\ccK$):

\begin{theo}
~~~~~~

\noindent 1. Let $G,H \in \ccG^{\circ}$.
\vspace{1 mm}

 a.  $G$ and $H$ have the same d-homotopy type
if, and only if, $C(G)$ and $C(H)$ have the same   homotopy type.

 b. $G$ and $H$ have the same d-homotopy type
if, and only if,  $\DG(G)$ and $\DG(H)$ have the same strong homotopy type.
\vspace{2 mm}

\noindent 2. Let $P,Q \in \ccP$.
\vspace{1 mm}

 a. $P$ and $Q$ have the same  homotopy type
if, and only if,  $Comp(P)$ and $Comp(Q)$ have the same d-homotopy type.

 b. $P$ and $Q$ have the same  homotopy type
if, and only if,  $\DP(P)$ and $\DP(Q)$ have the same  strong homotopy type.
\vspace{2 mm}

\noindent 3. Let $K,L \in \ccK$.
\vspace{1 mm}

 a. $K$ and $L$ have the same  strong homotopy type
if, and only if,  $\FG(K)$ and $\FG(L)$ have the same d-homotopy type.

 b. $K$ and $L$ have the same  strong homotopy type
if, and only if,  $\FP(K)$ and $\FP(L)$ have the same homotopy type.
\end{theo}

\proof All these equivalences are immediate corollaries of 
previous results: Theorems \ref{theo-red-C} and \ref{theo-C-red} (1.a),
Theorem \ref{theo-Delta}(1.b), 
Theorem \ref{theo-Comp} and Proposition \ref{wd=d} (2.a),
Theorem \ref{theo-DP}  and Proposition \ref{wd=d} (2.b),
Theorems \ref{theo-red-Gamma} and \ref{theo-Gamma-red} (3.a),
\cite[Theorem 4.14.b]{barmin09} and Theorem \ref{theo-FP}  (3.b).
\endproof

We recall that there is an operation of barycentric subdivision
either for graphs, for posets, or for simplicial complexes (\cite{bfj08})
verifying $Bd=C \circ Comp=\FG \circ \DG$ in $\ccG^{\circ}$,
$Bd=Comp\circ C=\FP\circ \DP$ in $\ccP$ and
$Bd=\DG \circ FG=\DP\circ \FP$ in $\ccK$.
\begin{prop} \label{passage-a-Bd}
~~

\noindent 1. Let $G,H \in \ccG$; then,
$G \redd H \Longleftrightarrow Bd(G) \redd Bd(H)$
\ss

\noindent 2. Let $K,L \in \ccK$; then,
$K \redred L \Longleftrightarrow Bd(K) \redred Bd(L)$
\ss

\noindent 3. Let $P,Q \in \ccP$; then,
$P \redwd Q \Longleftrightarrow Bd(P) \redd Bd(Q)$
\end{prop}

\proof The assertions 1 and 2 are corollaries of
Theorems \ref{theo-Delta}, \ref{theo-red-Gamma}
and \ref{theo-Gamma-red} by using, respectively,  $\FG\circ \DG=Bd$ (in $\ccG$)
and $\DG \circ \FG=Bd$ (in $\ccK$).
The assertion 3 is a consequence of
Theorems \ref{theo-Comp}, \ref{theo-red-C},
\ref{theo-C-red} and equality  $C\circ Comp=Bd$ (in $\ccP$).
\endproof

\begin{rmk}
If $L$ is  reduced to a vertex of $K$,
the assertion 2 of Proposition \ref{passage-a-Bd} 
is \cite[Theorem 4.15]{barmin09}.
 \end{rmk}

\section{Remarks about the \homm \:  complex}

Let $G,H\in \ccG$. The set of morphisms from $G$ to $H$
is the vertex set of the reflexive graph   $\homg(G,H)$
and is also the set of vertices (or 0-dimensional cells)
of the polyhedral complex  $\homm(G,H)$ 
(\cite{babsonkozlov},\cite{kozlov08})
 whose cells are indexed by  functions (which will be 
 called \textsl{indexing functions})
$\eta : V(G) \to 2^{V(H)}\setminus\{\emptyset\}$, 
such that if $(x,y)\in E(G)$,
then $\eta(x)\times \eta(y) \subset E(H)$. 

\begin{ex}\label{example-P3-K3}
We will illustrate the results of this section with the example given by
the path $G=P_3$ (i.e., $V(G)=\{0,1,2\}$ and  $0\sim 1\sim 2$)
and the complete graph $H=K_3$ (i.e., $V(K)=\{a,b,c\}$ and $a\sim b\sim c \sim a$).
The notation 
{\small \begin{tabular}{c|c}&$r$\\&$s$\\&$t$\end{tabular}} 
will indicate a morphism 
from $P_3$ to $K_3$ which sends 0 to $r$, 1 to $s$ and 2 to $t$.
There are 12 morphisms from $P_3$ to $K_3$:
\vspace{1 mm}

\noindent
\hspace{-7 mm}
{\small 
\begin{tabular}{cccccccccccc}
u&v&w&x&y&z&f&g&h&j&k&l\\ 
\begin{tabular}{c|c}&$a$\\&$c$\\&$a$\end{tabular}&
\begin{tabular}{c|c}&$b$\\&$c$\\&$b$\end{tabular}&
\begin{tabular}{c|c}&$b$\\&$a$\\&$b$\end{tabular}&
\begin{tabular}{c|c}&$c$\\&$a$\\&$c$\end{tabular}&
\begin{tabular}{c|c}&$c$\\&$b$\\&$c$\end{tabular}&
\begin{tabular}{c|c}&$a$\\&$b$\\&$a$\end{tabular}&
\begin{tabular}{c|c}&$a$\\&$c$\\&$b$\end{tabular}&
\begin{tabular}{c|c}&$b$\\&$c$\\&$a$\end{tabular}&
\begin{tabular}{c|c}&$b$\\&$a$\\&$c$\end{tabular}&
\begin{tabular}{c|c}&$c$\\&$a$\\&$b$\end{tabular}&
\begin{tabular}{c|c}&$a$\\&$b$\\&$c$\end{tabular}&
\begin{tabular}{c|c}&$c$\\&$b$\\&$a$\end{tabular}
\end{tabular}
}
\vspace{1 mm}

\noindent The graph $\homg(P_3,K_3)$ and the polyhedral complex $\homm(K_3,P_3)$
are represented in Figure \ref{homg-et-homm-P3-K3}.

\begin{figure}[h]
 \begin{center}
\psset{unit=0.16 cm}

\pspicture(-20,-16)(25,12)
\rput(12,0){$\bullet$}\rput(14,8){$\bullet$}\rput(6,10){$\bullet$}
\rput(4,2){$\bullet$}
\rput(-12,0){$\bullet$}\rput(-14,8){$\bullet$}\rput(-6,10){$\bullet$}
\rput(-4,2){$\bullet$}
\rput(0,-5){$\bullet$}\rput(0,-15){$\bullet$}
\rput(6,-10){$\bullet$}\rput(-6,-10){$\bullet$}
\psline(-6,10)(6,10)\psline(12,0)(6,-10)\psline(-12,0)(-6,-10)
\pspolygon(6,10)(14,8)(12,0)(4,2)
\pspolygon(-6,10)(-14,8)(-12,0)(-4,2)
\pspolygon(0,-5)(6,-10)(0,-15)(-6,-10)
\psline(6,10)(12,0)\psline(4,2)(14,8)
\psline(-6,10)(-12,0)\psline(-4,2)(-14,8)
\psline(-6,-10)(6,-10)\psline(0,-5)(0,-15)
\pscircle(6,10.7){0.7}\pscircle(14.5,8.5){0.7}\pscircle(12.5,-0.5){0.7}
\pscircle(3.5,1.5){0.7}
\pscircle(-6,10.7){0.7}\pscircle(-14.5,8.5){0.7}\pscircle(-12.5,-0.5){0.7}
\pscircle(-3.5,1.5){0.7}
\pscircle(0,-4.3){0.7}\pscircle(0,-15.7){0.7}
\pscircle(6.5,-10.5){0.7}\pscircle(-6.5,-10.5){0.7}
\uput[r](12.5,-0.5){$v$}\uput[ur](14.5,8.5){$g$}\uput[u](6,10.7){$u$}
\uput[dl](3.5,1.5){$f$}
\uput[l](-12.5,-0.5){$y$}\uput[ul](-14.5,8.5){$k$}\uput[u](-6,10.7){$z$}
\uput[dr](-3.5,1.5){$l$}
\uput[u](0,-4.3){$j$}\uput[d](0,-15.7){$h$}
\uput[r](6.5,-10.5){$w$}\uput[l](-6.5,-10.5){$x$}
\rput(-22,-10){$\homg(P_3,K_3)$}
\endpspicture
\pspicture(-15,-16)(20,12)
\rput(12,0){$\bullet$}\rput(14,8){$\bullet$}\rput(6,10){$\bullet$}
\rput(4,2){$\bullet$}
\rput(-12,0){$\bullet$}\rput(-14,8){$\bullet$}\rput(-6,10){$\bullet$}
\rput(-4,2){$\bullet$}
\rput(0,-5){$\bullet$}\rput(0,-15){$\bullet$}
\rput(6,-10){$\bullet$}\rput(-6,-10){$\bullet$}
\psline(-6,10)(6,10)\psline(12,0)(6,-10)\psline(-12,0)(-6,-10)
\pspolygon[fillstyle=solid,fillcolor=gray](6,10)(14,8)(12,0)(4,2)
\pspolygon[fillstyle=solid,fillcolor=gray](-6,10)(-14,8)(-12,0)(-4,2)
\pspolygon[fillstyle=solid,fillcolor=gray](0,-5)(6,-10)(0,-15)(-6,-10)
\uput[r](12,0){$v$}\uput[ur](14,8){$g$}\uput[u](6,10){$u$}
\uput[dl](4,2){$f$}
\uput[l](-12,0){$y$}\uput[ul](-14,8){$k$}\uput[u](-6,10){$z$}
\uput[dr](-4,2){$l$}
\uput[u](0,-5){$j$}\uput[d](0,-15){$h$}
\uput[r](6,-10){$w$}\uput[l](-6,-10){$x$}
\rput(22,-10){$\homm(P_3,K_3)$}
\endpspicture
\caption{The graph $\homg(P3,K3)$ and the polyhedral complex $\homm(P_3,K_3)$}
\label{homg-et-homm-P3-K3}
\end{center}
\end{figure}

\end{ex}

\subsection{$\homm(-,-)$ and $\homg(-,-)$}

For studying the 
polyhedral complex $\homm(G,H)$, it is usual to consider its face poset
$\FP(\homm(G,H))$ whose elements are all indexing functions 
with order given by
$\eta \leq \eta' $ if and only if 
$\eta(x)\subset \eta'(x)$ for all $x$ in $V(G)$.

Actually, there is a natural identification of $\FP(\homm(G,H))$ with a subposet
of $C(\homg(G,H))$, the poset of complete subgraphs of $\homg(G,H)$.
Indeed, let $\eta \in \FP(\homm(G,H))$
and for every vertex $x$ of $G$, let us choose an element 
$y_x \in \eta(x)$. Then the application $f:V(G)\to V(H),
x \mapsto y_x$ is actually a morphism from $G$ to $H$ ; 
such an application will be called an \textsl{associated morphism 
to $\eta$}. The set of all morphisms associated to $\eta$
will be called $\Psi(\eta)$. By definition of indexing functions,
 $\Psi(\eta)$ induces a complete subgraph of $\homg(G,H)$
and we get an injective poset map
$$\begin{array}{rcl}
\Psi~:~ \FP(\homm(G,H)) & \longrightarrow & C(\homg(G,H))\\
\eta & \mapsto & \Psi(\eta)
\end{array}
$$
which identifies $\FP(\homm(G,H))$ with a subposet of $C(\homg(G,H))$.

Now let $[f_1,f_2,\ldots,f_k]\in C(\homg(G,H))$
(i.e., the set $\{f_1,f_2,\ldots,f_k\}$ of morphisms 
from $G$ to $H$ induces a complete subgraph of $\homg(G,H)$).
We define the indexing function 
$\Phi([f_1,f_2,\ldots,f_k]) : 
V(G) \to 2^{V(H)}\setminus \{\emptyset\}$
 by 
$\Phi([f_1,f_2,\ldots,f_k])(x)=
\{f_1(x),f_2(x),\ldots,f_k(x)\}$ for all $x \in V(G)$.
This gives a morphism of posets:
$$\begin{array}{rcl}
\Phi~:~C(\homg(G,H)) & \hookrightarrow & \FP(\homm(G,H))\\
~ [f_1,f_2,\ldots,f_k] & \mapsto & \Phi([f_1,f_2,\ldots,f_k])
\end{array}
$$

\begin{prop}\label{homg-homm-poset}
Let $G,H \in \ccG$. By identifying $\FP(\homm(G,H))$ 
with a subposet of $C(\homg(G,H))$, we have:
\vspace{1 mm}

\centerline{$C(\homg(G,H)) \redd \FP(\homm(G,H))$}
\end{prop}

\proof First, we note that $\Phi \circ \Psi=1_{\FP(\homm(G,H))}$.
This implies that $(\Psi \circ \Phi)^2=\Psi \circ \Phi$,
i.e. $\Psi \circ \Phi : C(\homm(G,H))\to C(\homm(G,H))$ 
is a retraction
on $\FP(\homm(G,H))$ (identified with $\Psi(\homm(G,H))$).
Next, for all $[f_1,f_2,\ldots,f_k]\in C(\homg(G,H))$,
$[f_1,f_2,\ldots ,f_k] \subset \Psi \circ \Phi([f_1,f_2,\ldots ,f_k])$,
i.e. $1_{C(\homg(G,H))}\leq \Psi \circ \Phi$.
The conclusion follows from Proposition \ref{Pdemont-sur-Fix}.
\endproof

\begin{ex}
The posets obtained when $G=P_3$ and $K=K_3$ are drawed in 
Figures \ref{C-homg-P3-K3} and \ref{FP-Hom-P3-K3}.
The dismantling sequence: 
$fuv,guv,fgu,fgv,fg,uv,hwx,jwx,hjw,hjx,wx,hj,kyz,lyz,kly,klz,yz,kl$
illustrates the Proposition \ref{homg-homm-poset}.

\begin{figure}[h]
 \begin{center}
\psset{yunit=0.2 cm,xunit=0.3,nodesep=-3pt}

\pspicture(-25,0)(30,26)
\cnodeput[fillstyle=solid,fillcolor=white,linecolor=white](-17,24){fguv}{\small fguv}
\cnodeput[fillstyle=solid,fillcolor=white,linecolor=white](0,24){hjwx}{\small hjwx}
\cnodeput[fillstyle=solid,fillcolor=white,linecolor=white](16,24){klyz}{\small klyz}

\cnodeput[fillstyle=solid,fillcolor=white,linecolor=white](-21,16){fuv}{\tiny fuv}
\cnodeput[fillstyle=solid,fillcolor=white,linecolor=white](-18,16){guv}{\tiny guv}
\cnodeput[fillstyle=solid,fillcolor=white,linecolor=white](-15,16){fgu}{\tiny fgu}
\cnodeput[fillstyle=solid,fillcolor=white,linecolor=white](-12,16){fgv}{\tiny fgv}

\cnodeput[fillstyle=solid,fillcolor=white,linecolor=white](-4,16){hwx}{\tiny hwx}
\cnodeput[fillstyle=solid,fillcolor=white,linecolor=white](-1,16){jwx}{\tiny jwx}
\cnodeput[fillstyle=solid,fillcolor=white,linecolor=white](2,16){hjw}{\tiny hjw}
\cnodeput[fillstyle=solid,fillcolor=white,linecolor=white](5,16){hjx}{\tiny hjx}

\cnodeput[fillstyle=solid,fillcolor=white,linecolor=white](20,16){klz}{\tiny klz}
\cnodeput[fillstyle=solid,fillcolor=white,linecolor=white](17,16){kly}{\tiny kly}
\cnodeput[fillstyle=solid,fillcolor=white,linecolor=white](14,16){lyz}{\tiny lyz}
\cnodeput[fillstyle=solid,fillcolor=white,linecolor=white](11,16){kyz}{\tiny kyz}

\ncline{fguv}{fgu}\ncline{fguv}{fuv}\ncline{fguv}{fgv}\ncline{fguv}{guv}
\ncline{hjwx}{jwx}\ncline{hjwx}{hwx}\ncline{hjwx}{hjx}\ncline{hjwx}{hjw}
\ncline{klyz}{lyz}\ncline{klyz}{kyz}\ncline{klyz}{klz}\ncline{klyz}{kly}

\cnodeput[fillstyle=solid,fillcolor=white,linecolor=white](-22,8){uv}{\tiny uv}
\cnodeput[fillstyle=solid,fillcolor=white,linecolor=white](-19,8){fu}{\tiny fu}
\cnodeput[fillstyle=solid,fillcolor=white,linecolor=white](-17,8){gu}{\tiny gu}
\cnodeput[fillstyle=solid,fillcolor=white,linecolor=white](-15,8){fg}{\tiny fg}
\cnodeput[fillstyle=solid,fillcolor=white,linecolor=white](-13,8){fv}{\tiny fv}
\cnodeput[fillstyle=solid,fillcolor=white,linecolor=white](-11,8){gv}{\tiny gv}

\cnodeput[fillstyle=solid,fillcolor=white,linecolor=white](-8,8){vw}{\tiny vw}
\cnodeput[fillstyle=solid,fillcolor=white,linecolor=white](-5,8){wx}{\tiny wx}
\cnodeput[fillstyle=solid,fillcolor=white,linecolor=white](-3,8){hw}{\tiny hw}
\cnodeput[fillstyle=solid,fillcolor=white,linecolor=white](-1,8){jw}{\tiny jw}
\cnodeput[fillstyle=solid,fillcolor=white,linecolor=white](1,8){hj}{\tiny hj}
\cnodeput[fillstyle=solid,fillcolor=white,linecolor=white](3,8){hx}{\tiny hx}
\cnodeput[fillstyle=solid,fillcolor=white,linecolor=white](5,8){jx}{\tiny jx}
\cnodeput[fillstyle=solid,fillcolor=white,linecolor=white](8,8){xy}{\tiny xy}

\cnodeput[fillstyle=solid,fillcolor=white,linecolor=white](21,8){lz}{\tiny lz}
\cnodeput[fillstyle=solid,fillcolor=white,linecolor=white](19,8){kz}{\tiny kz}
\cnodeput[fillstyle=solid,fillcolor=white,linecolor=white](17,8){kl}{\tiny kl}
\cnodeput[fillstyle=solid,fillcolor=white,linecolor=white](15,8){ly}{\tiny ly}
\cnodeput[fillstyle=solid,fillcolor=white,linecolor=white](13,8){ky}{\tiny ky}
\cnodeput[fillstyle=solid,fillcolor=white,linecolor=white](11,8){yz}{\tiny yz}

\cnodeput[fillstyle=solid,fillcolor=white,linecolor=white](24,8){uz}{\tiny uz}

\ncline{fuv}{fu}\ncline{fuv}{uv}\ncline{fuv}{fv}
\ncline{guv}{gu}\ncline{guv}{uv}\ncline{guv}{gv}
\ncline{fgu}{fg}\ncline{fgu}{gu}\ncline{fgu}{fu}
\ncline{fgv}{gv}\ncline{fgv}{fv}\ncline{fgv}{fg}

\ncline{hwx}{wx}\ncline{hwx}{hx}\ncline{hwx}{hw}
\ncline{jwx}{wx}\ncline{jwx}{jx}\ncline{jwx}{jw}
\ncline{hjw}{jw}\ncline{hjw}{hw}\ncline{hjw}{hj}
\ncline{hjx}{jx}\ncline{hjx}{hx}\ncline{hjx}{hj}

\ncline{klz}{lz}\ncline{klz}{kz}\ncline{klz}{kl}
\ncline{kly}{ly}\ncline{kly}{ky}\ncline{kly}{kl}
\ncline{lyz}{yz}\ncline{lyz}{lz}\ncline{lyz}{ly}
\ncline{kyz}{yz}\ncline{kyz}{kz}\ncline{kyz}{ky}

\cnodeput[fillstyle=solid,fillcolor=white,linecolor=white](-22,0){u}{\small u}
\cnodeput[fillstyle=solid,fillcolor=white,linecolor=white](-18,0){f}{\small f}
\cnodeput[fillstyle=solid,fillcolor=white,linecolor=white](-14,0){g}{\small g}
\cnodeput[fillstyle=solid,fillcolor=white,linecolor=white](-10,0){v}{\small v}
\cnodeput[fillstyle=solid,fillcolor=white,linecolor=white](-6,0){w}{\small w}
\cnodeput[fillstyle=solid,fillcolor=white,linecolor=white](-2,0){h}{\small h}
\cnodeput[fillstyle=solid,fillcolor=white,linecolor=white](22,0){z}{\small z}
\cnodeput[fillstyle=solid,fillcolor=white,linecolor=white](18,0){l}{\small l}
\cnodeput[fillstyle=solid,fillcolor=white,linecolor=white](14,0){k}{\small k}
\cnodeput[fillstyle=solid,fillcolor=white,linecolor=white](10,0){y}{\small y}
\cnodeput[fillstyle=solid,fillcolor=white,linecolor=white](6,0){x}{\small x}
\cnodeput[fillstyle=solid,fillcolor=white,linecolor=white](2,0){j}{\small j}

\ncline{fu}{f}\ncline{fu}{u}
\ncline{uv}{v}\ncline{uv}{u}
\ncline{fv}{f}\ncline{fv}{v}
\ncline{gu}{u}\ncline{gu}{g}
\ncline{gv}{v}\ncline{gv}{g}
\ncline{fg}{f}\ncline{fg}{g}

\ncline{hw}{w}\ncline{hw}{h}
\ncline{wx}{x}\ncline{wx}{w}
\ncline{hx}{x}\ncline{hx}{h}
\ncline{jw}{w}\ncline{jw}{j}
\ncline{jx}{j}\ncline{jx}{x}
\ncline{hj}{j}\ncline{hj}{h}

\ncline{kl}{l}\ncline{kl}{k}
\ncline{kz}{k}\ncline{kz}{z}
\ncline{lz}{l}\ncline{lz}{z}
\ncline{ly}{y}\ncline{ly}{l}
\ncline{ky}{k}\ncline{ky}{y}
\ncline{yz}{z}\ncline{yz}{y}

\ncline{vw}{v}\ncline{vw}{w}
\ncline{xy}{x}\ncline{xy}{y}
\ncline{uz}{z}\ncline{uz}{u}

\endpspicture
\caption{$C(\homg(P3,K3))$}
\label{C-homg-P3-K3}
\end{center}
\end{figure}

\begin{figure}[h]
 \begin{center}
\psset{yunit=0.2 cm,xunit=0.3,nodesep=-3pt}

\pspicture(-25,-1)(30,26)

\cnodeput[fillstyle=solid,fillcolor=white,linecolor=white](-15,16){fguv}{$\bullet$}
\rput(-15,20){\small \begin{tabular}{|c}$a,b$\\$c$\\$a,b$\end{tabular}}
\cnodeput[fillstyle=solid,fillcolor=white,linecolor=white](2,16){hjwx}{$\bullet$}
\rput(2,20){\small \begin{tabular}{|c}$b,c$\\$a$\\$b,c$\end{tabular}}
\cnodeput[fillstyle=solid,fillcolor=white,linecolor=white](14,16){klyz}{$\bullet$}
\rput(14,20){\small \begin{tabular}{|c}$a,c$\\$b$\\$a,c$\end{tabular}}

\cnodeput[fillstyle=solid,fillcolor=white,linecolor=white](-19,8){fu}
{$\bullet${\small 1}}
\cnodeput[fillstyle=solid,fillcolor=white,linecolor=white](-17,8){gu}
{$\bullet${\small 2}}
\cnodeput[fillstyle=solid,fillcolor=white,linecolor=white](-13,8){fv}
{$\bullet${\small 3}}
\cnodeput[fillstyle=solid,fillcolor=white,linecolor=white](-11,8){gv}
{$\bullet${\small 4}}

\rput(-19,11){\tiny \begin{tabular}{|c}$a$\\$c$\\$a,b$\end{tabular}}

\cnodeput[fillstyle=solid,fillcolor=white,linecolor=white](-8,8){vw}{$\bullet$}
\cnodeput[fillstyle=solid,fillcolor=white,linecolor=white](-3,8){hw}
{$\bullet${\small 5}}
\cnodeput[fillstyle=solid,fillcolor=white,linecolor=white](-1,8){jw}
{$\bullet${\small 6}}
\cnodeput[fillstyle=solid,fillcolor=white,linecolor=white](3,8){hx}
{$\bullet${\small 7}}
\cnodeput[fillstyle=solid,fillcolor=white,linecolor=white](5,8){jx}
{$\bullet${\small 8}}
\cnodeput[fillstyle=solid,fillcolor=white,linecolor=white](8,8){xy}{$\bullet$}

\rput(-3,11){\tiny \begin{tabular}{|c}$b$\\$a$\\$b,c$\end{tabular}}
\rput(-8,11){\tiny \begin{tabular}{|c}$b$\\$a,c$\\$b$\end{tabular}}
\rput(8,11){\tiny \begin{tabular}{|c}$c$\\$a,b$\\$c$\end{tabular}}

\cnodeput[fillstyle=solid,fillcolor=white,linecolor=white](21,8){lz}
{$\bullet${\small 12}}
\cnodeput[fillstyle=solid,fillcolor=white,linecolor=white](18,8){kz}
{$\bullet${\small 11}}
\cnodeput[fillstyle=solid,fillcolor=white,linecolor=white](15,8){ly}
{$\bullet${\small 10}}
\cnodeput[fillstyle=solid,fillcolor=white,linecolor=white](12,8){ky}
{$\bullet${\small 9}}

\rput(12,11){\tiny \begin{tabular}{|c}$a,c$\\$b$\\$c$\end{tabular}}

\cnodeput[fillstyle=solid,fillcolor=white,linecolor=white](24,8){uz}{$\bullet$}
\rput(24,11){\tiny \begin{tabular}{|c}$a$\\$b,c$\\$a$\end{tabular}}

\ncline{fguv}{gu}\ncline{fguv}{fu}\ncline{fguv}{gv}\ncline{fguv}{fv}
\ncline{hjwx}{hw}\ncline{hjwx}{jw}\ncline{hjwx}{hx}\ncline{hjwx}{jx}
\ncline{klyz}{ky}\ncline{klyz}{ly}\ncline{klyz}{kz} \ncline{klyz}{lz}

\cnodeput[fillstyle=solid,fillcolor=white,linecolor=white](-22,0){u}{$\bullet$}
\cnodeput[fillstyle=solid,fillcolor=white,linecolor=white](-18,0){f}{$\bullet$}
\cnodeput[fillstyle=solid,fillcolor=white,linecolor=white](-14,0){g}{$\bullet$}
\cnodeput[fillstyle=solid,fillcolor=white,linecolor=white](-10,0){v}{$\bullet$}
\cnodeput[fillstyle=solid,fillcolor=white,linecolor=white](-6,0){w}{$\bullet$}
\cnodeput[fillstyle=solid,fillcolor=white,linecolor=white](-2,0){h}{$\bullet$}
\cnodeput[fillstyle=solid,fillcolor=white,linecolor=white](22,0){z}{$\bullet$}
\cnodeput[fillstyle=solid,fillcolor=white,linecolor=white](18,0){l}{$\bullet$}
\cnodeput[fillstyle=solid,fillcolor=white,linecolor=white](14,0){k}{$\bullet$}
\cnodeput[fillstyle=solid,fillcolor=white,linecolor=white](10,0){y}{$\bullet$}
\cnodeput[fillstyle=solid,fillcolor=white,linecolor=white](6,0){x}{$\bullet$}
\cnodeput[fillstyle=solid,fillcolor=white,linecolor=white](2,0){j}{$\bullet$}

\uput[d](-22,0){u}
\uput[d](-18,0){f}
\uput[d](-14,0){g}
\uput[d](-10,0){v}
\uput[d](-6,0){w}
\uput[d](-2,0){h}
\uput[d](22,0){z}
\uput[d](18,0){l}
\uput[d](14,0){k}
\uput[d](10,0){y}
\uput[d](6,0){x}
\uput[d](2,0){j}

\ncline{fu}{f}\ncline{fu}{u}
\ncline{fv}{f}\ncline{fv}{v}
\ncline{gu}{u}\ncline{gu}{g}
\ncline{gv}{v}\ncline{gv}{g}

\ncline{hw}{w}\ncline{hw}{h}
\ncline{hx}{x}\ncline{hx}{h}
\ncline{jx}{j}\ncline{jx}{x}
\ncline{jw}{j}\ncline{jw}{w}

\ncline{kz}{k}\ncline{kz}{z}
\ncline{lz}{l}\ncline{lz}{z}
\ncline{ly}{y}\ncline{ly}{l}
\ncline{ky}{k}\ncline{ky}{y}

\ncline{vw}{v}\ncline{vw}{w}
\ncline{xy}{x}\ncline{xy}{y}
\ncline{uz}{z}\ncline{uz}{u}

\endpspicture
\caption{$\FP(\homm(P3,K3))$}
\label{FP-Hom-P3-K3}
\end{center}
\end{figure}

\end{ex}

The \textsl{face graph}
$\FG(\homm(G,H)$ of the polyhedral complex $\homm(G,H)$ 
is the graph whose vertices  are the indexing functions  of $\homm(G,H)$
with edges $\eta \sim \eta'$ if and only if 
either $\eta(x)\subset \eta'(x)$ for all $x$ in $V(G)$,
or $\eta'(x)\subset \eta(x)$ for all $x$ in $V(G)$.
In other words, $\FG(\homm(G,H))=Comp(\FP(\homm(G,H)))$.

\begin{cor}\label{homg-homm-graph}
Let $G,H \in \ccG$. we have the following 
dismantling in $\ccG$ :
\vspace{1 mm}

\centerline{$Bd(\homg(G,H)) \redd \FG(\homm(G,H))$}
\end{cor}

\proof Follows from Proposition \ref{homg-homm-poset} by using
$Comp\circ C = Bd$, $Comp \circ \FP=\FG$ and Proposition \ref{dismant-PG}.
\endproof

\subsection{$\homm(G,H)$ and foldings in $G$ or in $H$}

\begin{theo}\label{homm-an-foldings-graph}
Let $G,H \in \ccG$.
\vspace{1 mm}

1. If $a$ is dismantlable in $G$, then 
$\FG(\homm(G,H))\redd \FG(\homm(G-a,H))$

(by identifying $\FG(\homm(G-a,H))$ 
with a subgraph of $\FG(\homm(G,H))$).
\vspace{2 mm}

2. If $u$ is dismantlable in $H$, 
then $\FG(\homm(G,H))\redd \FG(\homm(G,H-u))$

(by identifying $\FG(\homm(G,H-u))$ 
with a subgraph of $\FG(\homm(G,H))$).
\end{theo}

\proof
1. We have the following diagram where the morphisms
$A$ and $A'$ are dismantlings given by Corollary \ref{homg-homm-graph}
and the morphism B is a dismantling given by Propositions \ref{foldings-on-homg}
and \ref{passage-a-Bd}.1:
$$
\begin{psmatrix}[colsep=2cm,rowsep=2cm]
Bd(\homg(G,H)) & \FG(\homm(G,H))\\
Bd(\homg(G-a,H)) & \FG(\homm(G-a,H))
\psset{arrows=->,labelsep=5pt,nodesep= 0.5cm}
\ncline{1,1}{1,2}_{{\rm \bf A}}
\ncline{1,1}{2,1}<{{\rm \bf B}}
\ncline{2,1}{2,2}^{{\rm \bf A'}}
\ncline[linestyle=dashed]{1,2}{2,2}
\end{psmatrix}
$$
By considering   $(G,G',G'')=(Bd(\homg(G,H)),\FG(\homm(G,H)),\FG(\homm(G-a,H)))$,
the conclusion  follows from   Corollary \ref{dismant-subgraph}.

2. The proof is similar.
\endproof

\begin{ex}
Returning to the case $G=P_3$ and $H=K_3$, we have $2\domd 0$ in $P_3$
and $P_3 \redd P_3-2=K_2$. The deletion in $\FG(\homm(P_3,K_3))$
of  the twelve numbered vertices in the order indicated in 
Figure  \ref{FP-Hom-P3-K3}  followed by the deletion of the vertices
$f$, $g$, $h$, $j$, $k$ and $l$ is a dismantling sequence
from $\FG(\homm(P_3,K_3))$ to $\FG(\homm(K_2,K_3))$ (identified
with a subgraph of $\FG(\homm(P_3,K_3))$). We note that 
$\FG(\homm(K_2,K_3))$ is a stiff graph.

\begin{figure}[h]
 \begin{center}
\psset{unit=0.15 cm}

\pspicture(20,-18)(40,18)
\rput(12,0){$\bullet$}\rput(14,8){$\bullet$}\rput(6,10){$\bullet$}
\rput(4,2){$\bullet$}
\rput(-12,0){$\bullet$}\rput(-14,8){$\bullet$}\rput(-6,10){$\bullet$}
\rput(-4,2){$\bullet$}
\rput(0,-5){$\bullet$}\rput(0,-15){$\bullet$}
\rput(6,-10){$\bullet$}\rput(-6,-10){$\bullet$}
\psline(-6,10)(6,10)\psline(12,0)(6,-10)\psline(-12,0)(-6,-10)
\pspolygon(6,10)(14,8)(12,0)(4,2)
\pspolygon(-6,10)(-14,8)(-12,0)(-4,2)
\pspolygon(0,-5)(6,-10)(0,-15)(-6,-10)
\psline(6,10)(12,0)\psline(4,2)(14,8)
\psline(-6,10)(-12,0)\psline(-4,2)(-14,8)
\psline(-6,-10)(6,-10)\psline(0,-5)(0,-15)
\rput(9,5){$\bullet$}\rput(-9,5){$\bullet$}\rput(0,-10){$\bullet$}
\rput(9,-5){$\bullet$}\rput(-9,-5){$\bullet$}\rput(0,10){$\bullet$}

\rput(10,9){$\bullet$}\rput(8,1){$\bullet$}\psline(10,9)(8,1)
\rput(-10,9){$\bullet$}\rput(-8,1){$\bullet$}\psline(-10,9)(-8,1)
\rput(13,4){$\bullet$}\rput(5,6){$\bullet$}\psline(13,4)(5,6)
\rput(-13,4){$\bullet$}\rput(-5,6){$\bullet$}\psline(-13,4)(-5,6)
\rput(3,-7.5){$\bullet$}\rput(-3,-12.5){$\bullet$}\psline(3,-7.5)(-3,-12.5)
\rput(-3,-7.5){$\bullet$}\rput(3,-12.5){$\bullet$}\psline(-3,-7.5)(3,-12.5)

\pscircle(6,10.7){0.7}\pscircle(14.5,8.5){0.7}\pscircle(12.5,-0.5){0.7}
\pscircle(3.5,1.5){0.7}
\pscircle(-6,10.7){0.7}\pscircle(-14.5,8.5){0.7}\pscircle(-12.5,-0.5){0.7}
\pscircle(-3.5,1.5){0.7}
\pscircle(0,-4.3){0.7}\pscircle(0,-15.7){0.7}
\pscircle(6.5,-10.5){0.7}\pscircle(-6.5,-10.5){0.7}

\pscircle(0,10.7){0.7}\pscircle(9.6,-5.4){0.7}\pscircle(-9.6,-5.4){0.7}

\pscircle(3.5,-13){0.7}\pscircle(-3.5,-13){0.7}
\pscircle(3.5,-7){0.7}\pscircle(-3.5,-7){0.7}

\uput[dr](3.5,-13){5}\uput[dl](-3.5,-13){7}
\uput[ur](3.5,-7){6}\uput[ul](-3.5,-7){8}

\pscircle(10.1,9.65){0.7}\pscircle(7.8,0.35){0.7}
\pscircle(13.65,3.8){0.7}\pscircle(4.35,6.2){0.7}
\pscircle(-10.1,9.65){0.7}\pscircle(-7.8,0.35){0.7}
\pscircle(-13.65,3.8){0.7}\pscircle(-4.35,6.2){0.7}

\uput[u](10.1,9.65){2}\uput[d](7.8,0.35){3}
\uput[r](13.65,3.8){4}\uput[l](4.35,6.2){1}
\uput[u](-10.1,9.65){11}\uput[d](-7.8,0.35){10}
\uput[l](-13.65,3.8){9}\uput[r](-4.35,6.2){12}

\pscircle(9.6,5.4){0.7}
\pscircle(-9.6,5.4){0.7}
\pscircle(0,-10.7){0.7}
\pscircle(13.65,3.8){0.7}\pscircle(4.35,6.2){0.7}

\uput[r](12.5,-0.5){$v$}\uput[ur](14.5,8.5){$g$}\uput[u](6,10.7){$u$}
\uput[dl](3.5,1.5){$f$}
\uput[l](-12.5,-0.5){$y$}\uput[ul](-14.5,8.5){$k$}\uput[u](-6,10.7){$z$}
\uput[dr](-3.5,1.5){$l$}
\uput[u](0,-4.3){$j$}\uput[d](0,-15.7){$h$}
\uput[r](6.5,-10.5){$w$}\uput[l](-6.5,-10.5){$x$}
\endpspicture
\pspicture(0,-18)(15,18)
\rput(12,0){$\bullet$}\rput(6,10){$\bullet$}
\rput(-12,0){$\bullet$}\rput(-6,10){$\bullet$}
\rput(6,-10){$\bullet$}\rput(-6,-10){$\bullet$}
\psline(-6,10)(6,10)\psline(12,0)(6,-10)\psline(-12,0)(-6,-10)
\psline(6,10)(12,0)
\psline(-6,10)(-12,0)
\psline(-6,-10)(6,-10)
\rput(9,5){$\bullet$}\rput(-9,5){$\bullet$}\rput(0,-10){$\bullet$}
\rput(9,-5){$\bullet$}\rput(-9,-5){$\bullet$}\rput(0,10){$\bullet$}
\pscircle(6,10.7){0.7}\pscircle(12.5,-0.5){0.7}
\pscircle(-6,10.7){0.7}\pscircle(-12.5,-0.5){0.7}
\pscircle(6.5,-10.5){0.7}\pscircle(-6.5,-10.5){0.7}

\pscircle(0,10.7){0.7}\pscircle(9.6,-5.4){0.7}\pscircle(-9.6,-5.4){0.7}

\pscircle(9.6,5.4){0.7}
\pscircle(-9.6,5.4){0.7}
\pscircle(0,-10.7){0.7}
\rput(14,-2){$v'$}\uput[ur](6,10){$u'$}
\rput(-14,-2){$y'$}\uput[ul](-6,10){$z'$}
\rput(8.5,-10){$w'$}\rput(-8.5,-10){$x'$}
\rput(35,0){\small
\begin{tabular}{ccc}
$u'$ & $v'$ &$w'$ \\
\begin{tabular}{c|c}0&$a$\\1&$c$\end{tabular}&
\begin{tabular}{c|c}0&$b$\\1&$c$\end{tabular}&
\begin{tabular}{c|c}0&$b$\\1&$a$\end{tabular}\\ &&\\
$x'$ & $y'$ &$z'$\\
\begin{tabular}{c|c}0&$c$\\1&$a$\end{tabular}&
\begin{tabular}{c|c}0&$c$\\1&$b$\end{tabular}&
\begin{tabular}{c|c}0&$a$\\1&$b$\end{tabular}
\end{tabular}
}
\endpspicture
\caption{The graphs $\FG(\homm(P3,K3))$ and $\FG(\homm(K2,K3))$}
\label{FG-Hom-P3-K3-et-Hom-K2-K3}
\end{center}
\end{figure}
\end{ex}

Using $\DG \circ \FG=Bd$ (in $\ccK$) and Theorem \ref{theo-Delta}, 
we get the following corollary:

\begin{cor} Let $G,H \in \ccG$.
\vspace{1 mm}

\noindent 1. If $a$ is dismantlable in $G$, then 
$Bd(\homm(G,H))\redred Bd(\homm(G-a,H))$
\vspace{1 mm}

In particular, $Bd(\homm(G,H))$ and $Bd(\homm(G-a,H))$
have the same strong homotopy type.
\vspace{2 mm}

\noindent 2. If $u$ is dismantlable in $H$, then 
$Bd(\homm(G,H))\redred Bd(\homm(G,H-u))$
\vspace{1 mm}

In particular, $Bd(\homm(G,H)$ and $Bd(\homm(G,H-u))$
have the same strong homotopy type.
\end{cor}

\begin{rmk}
The result of the corollary implies the well known 
results about the behaviour of $\homm(G,H)$
in relation to foldings in $G$ or $H$ (\cite{babsonkozlov},
\cite{kozlov06a},\cite{csorba}).
It is actually more precise because the \textit{strong homotopy type} 
is very much stronger than the \textit{simple homotopy type}
(cf. Remark \ref{rigidity-d-homotopy})
and the equality of the two strong homotopy types follows from 
a strong collapse.
Moreover, this result in the framework of simplicial complexes
 is itself a consequence of Theorem \ref{homm-an-foldings-graph}
 in the framework of graphs.
\end{rmk}


\bibliographystyle{elsarticle-num}

\end{document}